\documentclass[reqno]{amsart}
\usepackage{mathtools,amssymb,wasysym,mathscinet,enumitem}
\usepackage{mathrsfs}  
\usepackage{geometry}
\numberwithin{equation}{section}

\usepackage[numbers,sort&compress]{natbib}
\usepackage{hyperref}
\usepackage{esint} 

\theoremstyle{plain}
\newtheorem{theorem}{Theorem}[section]
\newtheorem{lemma}[theorem]{Lemma}

\newtheorem{proposition}[theorem]{Proposition}

\theoremstyle{definition}
\newtheorem{definition}[theorem]{Definition}

\theoremstyle{remark}
\newtheorem{remark}[theorem]{Remark}

\makeatletter
\newcommand{\norm}{\@ifstar{\@normb}{\@normi}}
\newcommand{\@normb}[2]{\left\Vert{#1}\right\Vert_{#2}}
\newcommand{\@normi}[2]{\Vert{#1}\Vert_{#2}}
\newcommand{\norma}{\@ifstar{\@normba}{\@normia}}
\newcommand{\@normba}[2]{\left\vert{#1}\right\vert_{#2}}
\newcommand{\@normia}[2]{\vert{#1}\vert_{#2}}
\newcommand{\normb}{\@ifstar{\@normba}{\@normib}}
\newcommand{\@normbb}[2]{\left[{#1}\right]_{#2}}
\newcommand{\@normib}[2]{[{#1}]_{#2}}
\makeatother

\global\long\def\Sob#1#2{{W}^{#1}_{#2}} 
\global\long\def\dSob#1#2{{\mathcal{H}}^{#1}_{#2}} 
\global\long\def\tSob#1{{H}^{#1}}

\global\long\def\Sobloc#1#2{W^{#1}_{#2,\mathrm{loc}}} 
 
\global\long\def\oSob#1#2{\mathring{W}^{#1}_{#2}}
\global\long\def\Leb#1{L_{#1}} 
\global\long\def\Lebloc#1{L_{#1,\mathrm{loc}}}

\global\long\def\DMO{\mathrm{DMO}}
\newcommand{\action}[1]{\left<#1 \right>}

\DeclareMathOperator{\Div}{div}
\newcommand{\relphantom}[1]{\mathrel{\phantom{#1}}}

\newcommand{\myd}[1]{\,d{#1}}

\DeclareMathOperator{\supp}{supp}

\makeatletter
\def\@tocline#1#2#3#4#5#6#7{\relax
	\ifnum #1>\c@tocdepth % then omit
	\else
	\par \addpenalty\@secpenalty\addvspace{#2}%
	\begingroup \hyphenpenalty\@M
	\@ifempty{#4}{%
		\@tempdima\csname r@tocindent\number#1\endcsname\relax
	}{%
		\@tempdima#4\relax
	}%
	\parindent\z@ \leftskip#3\relax \advance\leftskip\@tempdima\relax
	\rightskip\@pnumwidth plus4em \parfillskip-\@pnumwidth
	#5\leavevmode\hskip-\@tempdima
	\ifcase #1
	\or\or \hskip 1em \or \hskip 2em \else \hskip 3em \fi%
	#6\nobreak\relax
	\hfill\hbox to\@pnumwidth{\@tocpagenum{#7}}\par% <---- \dotfill -> \hfill
	\nobreak
	\endgroup
	\fi}
\makeatother

\usepackage{tcolorbox}

\allowdisplaybreaks

\newcommand{\boldF}{\mathbf{F}}
\newcommand{\boldG}{\mathbf{G}}
\newcommand{\boldH}{\mathbf{H}}

\newcommand{\Rho}{\mathrm{P}}

\begin{document}
	
	\title[Stokes equations]{Spatial $C^1$, $C^2$, and Schauder estimates for nonstationary Stokes equations with Dini mean oscillation coefficients}
	
	\author[H. Dong]{Hongjie Dong}
	\address{Division of Applied Mathematics, Brown University, 182 George Street, Providence, RI 02912, USA}
	\email{hongjie\_dong@brown.edu }
	
	\author[H. Kwon]{Hyunwoo Kwon}
	\address{Division of Applied Mathematics, Brown University, 182 George Street, Providence, RI 02912, USA}
	\email{hyunwoo\_kwon@brown.edu }
	\thanks{H. Dong and H. Kwon were partially supported by the NSF under agreement DMS-2350129.}

	\begin{abstract}
		We establish the spatial differentiability of weak solutions to nonstationary Stokes equations in divergence form with variable viscosity coefficients having $L_2$-Dini mean oscillations. As a corollary, we derive local spatial Schauder estimates for such equations if the viscosity coefficient belongs to $C^\alpha_x$. Similar results also hold for strong solutions to nonstationary Stokes equations in nondivergence form.
	\end{abstract}

	\maketitle

	\section{Introduction}

	Consider the nonstationary Stokes equations with variable viscosity coefficients in divergence form: 
	\begin{equation}\label{eq:div-form}
		\left\{
		\begin{alignedat}{2}
			\partial_t u -D_i(a^{ij}D_j u)+\nabla \pi&=\Div \boldF&&\quad \text{in } U,\\
			\Div u&=g&&\quad \text{in } U.\\
		\end{alignedat}
		\right.
	\end{equation}
	Here $U$ {is} a cylindrical domain in $\mathbb{R}^{d+1}$, $u:U\rightarrow\mathbb{R}^d$ and $\pi:U\rightarrow\mathbb{R}$ denote the velocity field and the pressure of the fluid, respectively. Additionally, $\boldF:U\rightarrow\mathbb{R}^{d\times d}$ is a $2$-tensor, $\Div \boldF:U\rightarrow\mathbb{R}^d$ is a vector field defined by 
	\[ \Div \boldF =(D_j F^{1j},\dots, D_j F^{dj}), \]
	and $g:U\rightarrow\mathbb{R}$ is a function. The {coefficient matrix $A=(a^{ij})_{i,j=1}^d$}  satisfies the uniform ellipticity condition, i.e., there exists a constant $\nu \in (0,1)$ such that 
	\[  a^{ij}(t,x)\xi_i \xi_j\geq \nu |\xi|^2\quad \text{and}\quad |a^{ij}(t,x)|\leq \nu^{-1} \]
	hold for all $(t,x)\in\mathbb{R}^{d+1}$ and $\xi \in \mathbb{R}^d$. Throughout this paper, we adopt the Einstein summation convention for repeated indices unless stated otherwise.
	
	We also consider nonstationary Stokes equations in nondivergence form:
	\begin{equation}\label{eq:nondiv-form}
		\left\{
		\begin{alignedat}{2}
			\partial_t u -a^{ij}D_{ij} u+\nabla \pi&=f&&\quad \text{in } U,\\
			\Div u&=g&&\quad \text{in } U.
		\end{alignedat}
		\right.
	\end{equation}
	
Besides mathematical interests, our systems \eqref{eq:div-form} and \eqref{eq:nondiv-form} were motivated by the study of inhomogeneous fluids with density-dependent viscosity (see e.g. \cite{ABZ11}) and non-Newtonian fluids having time-dependent shear-thinning properties (see e.g. \cite{BGL16,BMR07}). Also, such systems arise in the study of Stokes equations on manifolds (see, e.g., \cite{DM04, MT01}). Recent studies have focused on local $L_{s, q}$-estimates for derivatives of solutions to these systems (see Dong-Kwon \cite{DK23, DK24} and references therein). This paper aims to provide a natural counterpart by establishing criteria that ensure spatial Schauder, $C^1$, and $C^2$ estimates for systems \eqref{eq:div-form} and \eqref{eq:nondiv-form}, thereby extending the regularity theory of parabolic equations to incorporate the nonlocal effects of the pressure.
	
	We briefly review the history of $C^1$-estimates for elliptic equations in divergence form to give a motivation to study our problem. Consider an elliptic equation whose {leading coefficients are} uniformly elliptic:
	\begin{equation}\label{eq:elliptic-div-form}
		-D_i(a^{ij}(x)D_j u)=D_i f^i\quad \text{in } B_1.
	\end{equation}
	
	Starting in the early 1930s, it is a classical result due to Schauder that the gradient of weak solutions to \eqref{eq:elliptic-div-form} belongs to $C^\alpha(B_{1/2})$ if {$A$} and $f$ are in $C^\alpha(B_1)$ for some $\alpha \in (0,1)$. Thus, it is natural to find a condition on the leading coefficients and external data for which the solution becomes $C^1$. In this context,  Li \cite{Li17} asked whether a weak solution to \eqref{eq:elliptic-div-form} is $C^1$ if $f=0$ and {
	\begin{equation}\label{eq:rho-function}  \rho(r):=\sup_{x \in B_{3/4}} \left(\sum_{i,j=1}^d\fint_{B_r(x)} |a^{ij}-(a^{ij})_{B_r(x)}|^2 \myd{y} \right)^{1/2},\quad 0<r\leq \frac{1}{4},
	\end{equation}
	satisfies the Dini condition, i.e., $\rho(0)=0$ and $\int_0^1 r^{-1} {\rho(r)}dr<\infty.$
	Here 
	\[(a^{ij})_{B_r(x)}=\fint_{B_r(x)} a^{ij}\myd{x} :=\frac{1}{|B_r|} \int_{B_r(x)} a^{ij}\myd{x}.\] This question was resolved by Dong and Kim \cite{DK17}, who replaced the $\Leb{2}$-norm in \eqref{eq:rho-function} with the $\Leb{1}$-norm.} Later, Dong, Escauriaza, and Kim \cite{DEK21} proved a similar result for parabolic equations. Throughout this paper, we will refer to such a condition as the Dini mean oscillation condition (see Definition \ref{defn:DMO}).
	
	For stationary Stokes equations in divergence form, Choi and Dong \cite{CD19b} established $C^1$ estimates and weak type-$(1,1)$ estimates for weak solutions under the condition that $a^{ij}$ is merely measurable in one direction, has Dini mean oscillation in the other direction, and $\boldF$ satisfies the Dini mean oscillation condition. Subsequently, they \cite{CD19} obtained global estimates on $C^{1,\mathrm{Dini}}$ domains. These results were further extended by Choi, Dong, and Xu \cite{CDX22} to cases where the viscosity coefficients have piecewise Dini mean oscillation. In addition, Dong, Li, and Xu \cite{DLX24} achieved higher regularity for weak solutions when the viscosity coefficients possess piecewise H\"older continuity.
	
	However, there are few corresponding regularity results on nonstationary Stokes equations even if $a^{ij}=\delta^{ij}$. Solonnikov \cite{Son00} first obtained global Schauder estimates for \eqref{eq:nondiv-form} in $(0,T)\times\mathbb{R}^3_+$. Later, Chang and Kang \cite{CK18b} proved the H\"older continuity of solutions to Stokes equations under suitable compatibility conditions on the boundary data. Related to our work, very recently, independent to our work, Dong, Li, and Wang \cite{DLW24} obtained local spatial Schauder estimates for \eqref{eq:nondiv-form} when {$A \in C^\alpha_x$}, $f\in C^\alpha_x(Q_1)$,  and $g\in C^{1,\alpha}_x(Q_1)$. More precisely, they proved that if $u\in \Sob{1,2}{2}(Q_1)$ is a strong solution of \eqref{eq:nondiv-form} in $Q_1$, then 
	\[ \norm{D^2u}{C^\alpha_x(Q_{1/4})}\leq N\left(\norm{u}{\Leb{\infty}\Leb{2}(Q_1)}+\norm{u}{\Sob{0,2}{2}(Q_1)}+\norm{f}{C^\alpha_x(Q_1)}+\norm{g}{C^{1,\alpha}_x(Q_1)}\right) \]
	for some constant $N=N(d,\nu,\alpha)>0$. See \eqref{eq14.00} for the definition of the $C^\alpha_x$ norm.

	The purpose of this paper is to obtain $C^1_x$-estimates for \eqref{eq:div-form} and $C^2_x$-estimates for \eqref{eq:nondiv-form} under suitable assumptions on variable viscosity coefficients and external data. For the equations in divergence form, we prove that the gradient of a weak solution to \eqref{eq:div-form} in $Q_1$ is bounded in $Q_{1/4}$ when $A$, $f$, and $g$ satisfy appropriate Dini mean oscillation conditions. Furthermore, $Du$ is continuous in $B_{1/8}$ for each $t \in (-1/64, 0)$, and the vorticity is continuous in both $t$ and $x$ {(Remark \ref{rem:continuity-vorticity})}. As a corollary, our method yields local spatial Schauder estimates (Remark \ref{rem:divergence-Holder}). A similar result holds for equations in nondivergence form. Precise statements of these results are provided in Theorems \ref{thm:A} and \ref{thm:B}, respectively.
	
	\subsubsection*{Outline of the proofs} 
	Let us outline the proofs of the main theorems. A natural approach is to use Campanato's method, as in \cite{G83, L96}. Broadly speaking, if we can demonstrate that the mean oscillation of $Du$ (or $D^2u$, respectively) in cylinders vanishes at a certain rate as the radii of the cylinders shrink to zero, then $Du$ becomes uniformly continuous in $(t,x)$, as shown in \cite{D12}. If the viscosity coefficients are smooth, one might expect $Du$ to be H\"older continuous. However, the following example by Serrin \cite{S62} illustrates that H\"older continuity in time cannot be achieved unless additional regularity is imposed on $u$. Specifically, consider the example:
	\begin{equation}\label{eq:Serrin} u(t,x) = c(t)\nabla h(x) \quad \text{and} \quad \pi(t,x) = -c'(t)h(x), \end{equation}
	where $c$ is an absolutely continuous function on $[-1,0]$ and $h$ is harmonic in $B_1$. In this case, $u$ is smooth in $x$ but not H\"older continuous in $t$. This demonstrates that the standard approach described above is insufficient.
	
	To address this issue, we first rewrite the original equation as
	\[ \partial_t u -D_i(\hat{a}^{ij}(t)D_j u)+\nabla \pi =\Div \boldF+D_i((a^{ij}-\hat{a}^{ij})D_ju)\quad \text{and}\quad \Div u =g\quad \text{in } Q_r,\]
	where $\hat{a}^{ij}(t)=\frac{1}{|B_r|}\int_{B_r} a^{ij}(t,x)\myd{x}$. Taking the curl to the above equation leads to a vorticity equation whose leading coefficients depend only on $t$. Then we obtain a decay rate of mean oscillation of the vorticity $\omega_{ij}=D_j u^i-D_i u^j$ via a standard argument; see e.g. \cite{D12}. Next, using the identity
	\begin{equation*}
		\Delta u^l= D_i\omega_{li} + D_l g,\quad l=1,\dots,d,
	\end{equation*}
	and Schauder estimates for the Poisson equation, we estimate the rate at which $u$ can be approximated by a polynomial while fixing $t$. This is achieved by analyzing the decay rate of the following quantity:
	\begin{equation}\label{eq:phi-rate-intro}
		\phi_1(r,t_0,x_0)= \frac{1}{r^{d/2+1}}\sup_{t\in (t_0-r^2,t_0)} \inf_{p \in \mathcal{P}_1} \norm{u(t,\cdot)-p}{\Leb{2}(B_r(x_0))},
	\end{equation}
	where $\mathcal{P}_1$ denotes the set of polynomials up to order $1$. The desired result then follows by an iterative argument, as in \cite{DK17, DEK21}. A similar approach is used to handle equations in nondivergence form 
	
	{A} closely related result appears in the independent work of Dong, Li, and Wang \cite{DLW24}, where the leading coefficients are assumed to belong to $C^\alpha_x$. In that proof, the rate of approximation of the velocity field by polynomials with \emph{time-dependent} coefficients was estimated to be of order $2+\alpha$ {(of order $1+\alpha$ if the equation is in divergence form)}. This required a meticulous selection of parameters to facilitate a delicate scaling and induction argument. The approach ultimately yielded the desired spatial Schauder estimates via a Campanato-type characterization. However, extending this method to the more general case where the leading coefficients satisfy only the Dini mean oscillation condition would present significant challenges. 
	
	In contrast to \cite{DLW24}, we adopt a more direct approach by analyzing the quantity \eqref{eq:phi-rate-intro}, which avoids the complex induction arguments employed in \cite{DLW24}. This alternative method accommodates the broader setting of coefficients with Dini mean oscillation and enables us to derive $C^1_x$ and $C^2_x$ estimates. Furthermore, as a corollary, it provides spatial Schauder estimates in this generalized framework.
	
	\subsubsection*{Organizations} 
	The rest of this paper is organized as follows. In Section \ref{sec:results}, we introduce notation and main results of this paper. Section \ref{sec:prelim} contains several estimates related to the Dini mean oscillation condition and Caccioppoli estimates for Stokes equations with simple coefficients. In Section \ref{sec:vorticity-rate}, we obtain the decay rate of the mean oscillation of the vorticity. In Section \ref{sec:velocity-rate}, we estimate the rate of the approximation of the velocity by polynomials, which will be crucially used in obtaining the main theorems. The proofs of Theorems \ref{thm:A} and \ref{thm:B} are given in Section \ref{sec:proof-main}. In Appendix \ref{sec:Appendix}, we prove the solvability of Stokes equations that we need. Finally, we present an approximation argument which will be used in the proofs of main theorems in Appendix \ref{app:approximation} for the sake of completeness.
	
	\section{Notation and main results}\label{sec:results}
	\subsection{Notation}
	%To explain main results, we introduce several concepts. 
	For $A\subset \mathbb{R}^d$ with $|A|<\infty$, where $|\cdot|$ denotes the $d$-dimensional Lebesgue measure, we write
	\[ \fint_A f \myd{x}:=\frac{1}{|A|} \int_A f \myd{x}. \]
	In particular, for $X_0=(t_0,x_0)\in\mathbb{R}\times\mathbb{R}^d$, we define
	\[ (f)_{X_0,r}=\fint_{Q_r(X_0)} f \myd{xdt} \quad \text{and}\quad  [f]_{x_0,r}(t)=\fint_{B_r(x_0)} f(t,x)dx.\]
	
	For a vector field $u=(u^1,\dots,u^d)$, we define the gradient and the vorticity of $u$ by 
	\[ (\nabla u)^{ij}= D_j u^i \quad \text{and}\quad \omega_{ij} =(\nabla \times u)^{ij}=D_j u^i -D_i u^j,\quad i,j=1,\dots,d,\]
	respectively. For a function $\phi$, we define 
	\[ (\nabla^2 \phi)^{ij}=D_{ij} \phi,\quad i,j=1,\dots,d.\]

	For two tensors $\boldF=(F^{ij})$ and $\boldG=(G^{ij})$,  their inner product is defined by
	\[ \boldF:\boldG = F^{ij} G^{ij}.\]

	By $N=N(p_1,\dots,p_k)$, we denote a generic positive constant depending only on the parameters $p_1$, ..., $p_k$. For $X_0=(t_0,x_0)\in \mathbb{R}^{d+1}$ and $r>0$, we write $B_r(x_0)$ the open ball centered at $x_0$ with radius $r>0$ and 
	\[ Q_r(X_0)=(t_0-r^2,t_0)\times B_r(x_0) \]
	the parabolic cylinder centered at $X_0$ with radius $r>0$. 
	
	By $C_0^\infty(U)$, we denote the space of infinitely differentiable functions with compact support in $U$. For $k\in\mathbb{N}$ and $1<q<\infty$, we write $\Sob{k}{q}(\Omega)$ the Sobolev space. When $q=2$, we write $\tSob{k}(\Omega)=\Sob{k}{2}(\Omega)$. We write $\oSob{1}{q}(\Omega)$ the closure of $C_0^\infty(\Omega)$ under $\Sob{1}{q}$-norm. 
	
	For $1<q<\infty$, we define parabolic Sobolev spaces
	\begin{align*}
		\Sob{0,1}{q}((S,T)\times\Omega)&:=\{u : u, Du\in \Leb{q}((S,T)\times\Omega)\},\\
		\Sob{1,2}{q}((S,T)\times\Omega)&:=\{u : u, Du, D^2 u, u_t\in \Leb{q}((S,T)\times\Omega)\}.
	\end{align*}
	We write $\Sob{0,1}{2,\sigma}((S,T)\times\Omega)$ the space of all vector fields in $\Sob{0,1}{2}((S,T)\times\Omega)$ with divergence free. 
	
	For equations in divergence form, we introduce another function spaces $\mathbb{H}^{-1}_q$ and $\dSob{1}{q}$. We say that $f \in \mathbb{H}^{-1}_q((S,T)\times\Omega)$ if there exist $g_0,g=(g_1,\dots,g_d) \in \Leb{q}((S,T)\times\Omega)$ such that 
	\[ f= g_0+ D_i g_i\quad \text{in } (S,T)\times \Omega \]
	in the sense of distributions and the norm 
	\[ 
	\norm{f}{\mathbb{H}^{-1}_q((S,T)\times\Omega)} :=\inf\left\{ \sum_{i=0}^d \norm{g_i}{\Leb{q}((S,T)\times\Omega)} : f = g_0 + D_i g_i \right\}
	\]
	is finite. We define 
	\[ \dSob{1}{q}((S,T)\times\Omega):=\{ u : u\in \Sob{0,1}{q}((S,T)\times\Omega), u_t \in \mathbb{H}^{-1}_q((S,T)\times\Omega)\} \]
	with the norm 
	\[ \norm{u}{\mathcal{H}^1_q((S,T)\times\Omega)} := \norm{u_t}{\mathbb{H}^{-1}_q((S,T)\times\Omega)} + \norm{u}{\Sob{0,1}{q}((S,T)\times\Omega)}.
	\]
	
	For $k\in\mathbb{N}$ and $\alpha,\beta \in (0,1]$, we define the H\"older semi-norm
	\[ [u]_{C^{0,\alpha}(\Omega)} = \sup_{x\neq y; x,y\in \Omega} \frac{|u(x)-u(y)|}{|x-y|^\alpha} \]
	and the H\"older norm 
	\[ \norm{u}{C^{k,\alpha}(\Omega)}=\sum_{|\gamma|\leq k}\norm{D^\gamma u}{\Leb{\infty}(\Omega)}+\sup_{|\gamma|=k}[D^\gamma u]_{C^{0,\alpha}(\Omega)}.
	\]
	
	Similarly, for $U=(S, T)\times \Omega$, we define the anisotropic H\"older semi-norm
	\[ [u]_{C^{\beta,\alpha}_{t,x}(U)} = \sup_{(t,x),(s,y)\in U; (t,x)\neq(s, y)} \frac{|u(t,x)-u(s,y)|}{|t-s|^{\beta}+|x-y|^\alpha} \]
	and the corresponding H\"older norm
	\[ \norm{u}{C^{\beta,\alpha}_{t,x}(U)}=\norm{u}{\Leb{\infty}(U)}+[u]_{C^{\beta,\alpha}(U)}.
	\]
	We also define 
	\[
	[u]_{C^\alpha_x(U)}=\sup_{(t,x),(t,y)\in U;x\neq y} \frac{|u(t,x)-u(t,y)|}{|x-y|^\alpha}
	\]
	and
	\begin{equation}
		\label{eq14.00}
		\norm{u}{C^\alpha_x(U)}=\norm{u}{\Leb{\infty}(U)}+[u]_{C^\alpha_x(U)}
	\end{equation}
	Similarly, we can define $C^{1,\alpha}_x(U)$.
	
	For $(S,T)\times \Omega \subset \mathbb{R}\times \mathbb{R}^d$, $k\in\mathbb{N}$, and $1\leq s,q\leq \infty$, we write
	\[ \Leb{s}\Leb{q}((S,T)\times \Omega)=\Leb{s}(S,T;\Leb{q}(\Omega)),\quad \Leb{s}\Sob{k}{q}((S,T)\times\Omega)=\Leb{s}(S,T;\Sob{k}{q}(\Omega)).\]
	When $s=q$, we write $\Leb{s}((S,T)\times\Omega)=\Leb{s}\Leb{s}((S,T)\times\Omega)$.

	Finally, we write $A\apprle B$ if there exists a constant $N$ independent of $A$ and $B$ such that $A\leq N B$. {We also write $A \apprle_{\alpha,\beta,\gamma...} B$ if the constant $N$ depends on the parameters $\alpha$, $\beta$, $\gamma$, ...}.
	
	\subsection{Main results}
	We present the main results of this paper in this subsection. To begin with, we introduce several definitions which will be used in this paper.

	\begin{definition}\label{defn:DMO} Let $f ,g\in \Leb{2}(Q_2)$. 
		\begin{enumerate}[label=\textnormal{(\roman*)}]
			\item We say that $f$  satisfies the $\Leb{2}$-\emph{Dini mean oscillation $(\DMO_x)$ in small cylinders} if
			\begin{align*}
				\rho_f(r)&=\sup_{(t_0,x_0)\in Q_1} \left(\fint_{Q_r(t_0,x_0)} |f-[f]_{x_0,r}(t)|^2 \myd{x}dt\right)^{1/2},\quad 0<r\leq 1,
			\end{align*}
			satisfies 
			\[ I_{\rho_f}(t)=\int_0^t\frac{\rho_f(s)}{s}ds<\infty\quad \text{for each } t\in (0,1].\]
			\item We say that $g$ satisfies the $\Leb{2}$-$\DMO_x$ \emph{in small balls} if
			\[
			\widehat{\rho}_g(r)=\sup_{(t,x_0)\in Q_1}  \left(\fint_{B_r(x_0)} |g(t,x)-[g]_{x_0,r}(t)|^2 \myd{x}\right)^{1/2},\quad 0<r\leq 1,
			\]
			satisfies $I_{\widehat{\rho}_g}(t)<\infty$ for each $t\in (0,1]$.
		\end{enumerate}
	\end{definition}

	We next give the definition of solutions that we consider in this paper.
	\begin{definition}
		For $-\infty\leq S<T\leq \infty$, let $U=(S,T)\times \Omega$, where $\Omega \subset \mathbb{R}^d$.
		\begin{enumerate}[label=\textnormal{(\roman*)}]
			\item We say that $u\in \Sobloc{0,1}{1}(U)$ is a weak solution to \eqref{eq:div-form} in $U$ if $u$ satisfies
			\[ \int_U u\cdot (-\partial_t \phi)\myd{x}dt+\int_U a^{ij}D_j u \cdot  D_i \phi \myd{x}dt=-\int_U \boldF:\nabla \phi \myd{x}dt \]
			for all $\phi \in C_0^\infty(U)$ with $\Div \phi =0$ in $U$ and 
			\[  -\int_\Omega u(t,x)\cdot \nabla \psi (x)\myd{x}=\int_\Omega g(t,x) \psi(x) \myd{x} \]
			for a.e. $t\in (S,T)$ and for all $\psi \in C_0^\infty(\Omega)$.
			\item We say that $(u,\pi)\in \Sobloc{1,2}{1}(U)\times \Sobloc{0,1}{1}(U)$ is a strong solution {to \eqref{eq:nondiv-form} in $U$} if $(u,\pi)$ satisfies \eqref{eq:nondiv-form} in $U$ pointwise a.e.
		\end{enumerate}
	\end{definition}
	\begin{remark}
		By de Rham's theorem, if $u\in \Leb{\infty}\Leb{2}(Q_1)\cap\Leb{2}\tSob{1}(Q_1)$ is a weak solution of \eqref{eq:div-form} in $Q_1$, then there exists a distribution $\pi(t)\in \Lebloc{2}(B_1)$ for almost any $t\in (-1,0)$ satisfying \eqref{eq:div-form} in the sense of distribution (see e.g. \cite[Theorem 2.1]{G00}). 
	\end{remark}

	Now we are ready to present our main results. The first result concerns the spatial differentiability of weak solutions to \eqref{eq:div-form} when $A$, $\boldF$, and $g$ satisfy $\Leb{2}$-$\DMO_x$.
	
	\begin{theorem}\label{thm:A}
		Let $\boldF\in \Leb{2}(Q_2)$ and $g\in \Leb{2}(Q_2)$. Suppose that $A$ and $\boldF$ satisfy the $\Leb{2}$-$\DMO_x$ in small cylinders and $g$ satisfies $\Leb{2}$-$\DMO_x$ in small balls. If $u \in \Leb{\infty}\Leb{2}(Q_1)\cap \Leb{2}\tSob{1}(Q_1)$ is a weak solution to \eqref{eq:div-form} in $Q_1$, then $Du \in \Leb{\infty}(Q_{1/4})$ and it is uniformly continuous in $B_{1/8}$ for each $t\in (-1/64,0)$. More precisely, for $\beta \in (0,1)$, we have
		\begin{align*}
			\norm{Du}{\Leb{\infty}(Q_{1/4})}\apprle \norm{u}{\Leb{\infty}\Leb{2}(Q_1)}+\norm{\boldF}{\Leb{2}(Q_1)}+\norm{g}{\Leb{2}(Q_1)}+I_{\Rho_{\boldF}}(1)+I_{\widehat{\Rho}_{g}}(1)
		\end{align*}
		and 
		\begin{equation}\label{eq:modulus-continuity-Du}
			\begin{aligned}
				&|Du(t,x)-Du(t,y)|\\
				&\apprle |x-y|^\beta \left(\norm{u}{\Leb{\infty}\Leb{2}(Q_1)}+\norm{\boldF}{\Leb{2}(Q_1)}+\norm{g}{\Leb{2}(Q_1)}\right)\\
				&\relphantom{=}+I_{\Rho_{\boldF}}(|x-y|)+I_{\widehat{\Rho}_{g}}(|x-y|)\\
				&\relphantom{=}+\left(\norm{u}{\Leb{\infty}\Leb{2}(Q_1)}+\norm{\boldF}{\Leb{2}(Q_1)}+\norm{g}{\Leb{2}(Q_1)}+I_{\Rho_{\boldF}}(1)+I_{\widehat{\Rho}_g}(1) \right)I_{\Rho_A}(|x-y|),
			\end{aligned}
		\end{equation}
		for $t\in (-1/64,0)$ and $x,y \in B_{1/8}$ with $|x-y|<1/8$, 
		where the implicit constant depends on $d$, $\nu$, $\beta$, and $A$.
		Here $\Rho_{\boldF}$ and $\widehat{\Rho}_{g}$ are defined in \eqref{eq:Rho-defn} and \eqref{eq:Rho-hat-g-defn}, respectively.
	\end{theorem}

	The second result concerns spatial twice differentiability of strong solutions to \eqref{eq:nondiv-form} when $A$, $f$, and $Dg$ satisfy $\Leb{2}$-$\DMO_x$. To state this result, let us write
	\[ 
	\tilde{W}^{1,2}_2(Q_1)=\{ u : u,Du,D^2u\in\Leb{2}(Q_1), u_t\in\Leb{1}(Q_1)\}.
	\]
	
	\begin{theorem}\label{thm:B}
		Let $f\in \Leb{2}(Q_2)$ and $g\in \Sob{0,1}{2}(Q_2)$. Suppose that $A$ and $f$ satisfy $\Leb{2}$-$\DMO_x$ in small cylinders and $Dg$ satisfies $\Leb{2}$-$\DMO_x$ in small balls. If $(u,\pi)\in \tilde{W}^{1,2}_2(Q_1)\times \Sob{0,1}{1}(Q_1)$ is a strong solution {to \eqref{eq:nondiv-form} in $Q_1$,} then $D^2u \in \Leb{\infty}(Q_{1/4})$ and it is uniformly continuous in $B_{1/8}$ for each $t\in (-1/64,0)$. More precisely, for $\beta\in(0,1)$, we have
		\begin{equation*}
			\begin{aligned}
				\norm{D^2u}{\Leb{\infty}(Q_{1/4})}&\apprle \norm{u}{\Leb{\infty}\Leb{2}(Q_{1})}+\norm{f}{\Leb{2}(Q_1)}+\norm{Dg}{\Leb{2}(Q_1)}+I_{\Rho_f}(1)+I_{\widehat{\Rho}_{Dg}}(1)
			\end{aligned}
		\end{equation*}
		and 
		\begin{equation*}
			\begin{aligned}\label{eq:Du-uniformly-conti-nondiv}
				&|D^2 u(t,x)-D^2u(t,y)|\nonumber\\
				&\apprle |x-y|^\beta\left( \norm{u}{\Leb{\infty}\Leb{2}(Q_1)}+\norm{f}{\Leb{2}(Q_1)}+\norm{Dg}{\Leb{2}(Q_1)}\right)\\
				&\relphantom{=}+I_{\Rho_f}(|x-y|)+I_{\widehat{\Rho}_{Dg}}(|x-y|)\\
				&\relphantom{=}+\left(\norm{u}{\Leb{\infty}\Leb{2}(Q_1)}+\norm{f}{\Leb{2}(Q_1)}+\norm{Dg}{\Leb{2}(Q_1)}+I_{\Rho_f}(1)+I_{\widehat{\Rho}_{Dg}}(1) \right)I_{\Rho_A}(|x-y|)\nonumber
			\end{aligned}
		\end{equation*}
		for $t\in (-1/64,0)$ and $x,y \in B_{1/8}$ with $|x-y|<1/8$, 
		where the implicit constant depends on $d$, $\nu$, $\beta$, and $A$.
	\end{theorem}
	
	\begin{remark}\leavevmode
		\begin{enumerate}[label=\textnormal{(\roman*)}]
			\item If $A\in C^\alpha_x$, $\boldF,g\in C^\alpha_x(Q_2)$ for some $\alpha \in (0,1)$, then it follows from \eqref{eq:modulus-continuity-Du} that $Du \in C^\alpha_x(Q_{1/4})$ and
			\[
			[Du]_{C^{\alpha}_x(Q_{1/8})}\apprle \norm{u}{\Leb{\infty}\Leb{2}(Q_1)}+\norm{\boldF}{C^\alpha_x(Q_1)}+\norm{g}{C^\alpha_x(Q_1)}.\]
			Moreover, the vorticity $\omega$ of $u$ belongs to $C^{\alpha/2,\alpha}_{t,x}(Q_{1/8})$ and satisfies 
			\begin{align*}
				[\omega]_{C^{\alpha/2,\alpha}_{t,x}(Q_{1/8})}&\apprle \norm{u}{\Leb{\infty}\Leb{2}(Q_1)}+\norm{\boldF}{C^\alpha_x(Q_1)}+\norm{g}{C^\alpha_x(Q_1)};\end{align*}
			see Remark \ref{rem:divergence-Holder}. {A similar result also holds for the equation in nondivergence form; see Remark \ref{rem:nondivergence-Holder}. We note that these results were recently proved by Dong-Li-Wang \cite{DLW24}.}
			\item It is a classical result that Hessian of solutions to the second-order parabolic equation has H\"older regularity in $t$ as well even if coefficients and data are only H\"older continuous in $x$ (see e.g. \cite[p.208]{K08}). However, for the Stokes equations, we cannot get better regularity in $t$ as suggested in Serrin's example \eqref{eq:Serrin}.  
			\item Suppose in addition that $\boldF,g\in C^{\alpha/2,\alpha}_{t,x}(Q_2)$ and $A\in C^{\alpha/2,\alpha}_{t,x}$ for some $\alpha \in (0,1)$. If $u\in C^{\alpha/2,\alpha}_{t,x}(Q_1)$ is a solution of \eqref{eq:div-form} in $Q_1$, then $Du \in C^{\alpha^-/2,\alpha}_{t,x}(Q_{1/4})$; see Remark \ref{rem:further-regularity}. A similar result holds for the equation in nondivergence form; see Remark \ref{rem:nondivergence-Holder-additional}.
		\end{enumerate}
	\end{remark}

	\section{Preliminaries}\label{sec:prelim}

	\subsection{Basic estimates}
	We first summarize estimates related to Dini functions (see e.g. \cite[Section 8]{CD19b}, \cite[Lemma 1]{D12}, \cite[Lemma 2.7]{DK17}).
	\begin{lemma}\label{lem:Dini-estimates}
		Let $\rho:[0,1]\rightarrow[0,\infty)$ satisfy
		\[ I_\rho(1)=\int_0^1 \frac{\rho(r)}{r}dr<\infty.\]
		\begin{enumerate}[label=\textnormal{(\roman*)}]
			\item Define
			\begin{equation}\label{eq:rho-tilde}
				\tilde{\rho}(t):=\sum_{k=0}^\infty a^k (\rho(b^k t)[b^k t\leq 1] +\rho(1)[b^k t>1]) 
			\end{equation}
			for some constants $a\in (0,1)$ and $b>1$. 
			Here we used the Iverson bracket notation, i.e., $[P]=1$ if $P$ is true and $[P]=0$ otherwise. Then 
			\[ \int_0^1\frac{\tilde{\rho}(r)}{r}dr<\infty.\]
			\item Let $\rho$ be a nonnegative bounded function. Suppose that there is $c_1,c_2>0$ and $0<\kappa<1$ such that 
			\begin{equation}\label{eq:dini-control}
				c_1\rho(t)\leq \rho(s)\leq c_2 \rho(t)\quad \text{whenever } \kappa t\leq s\leq t\quad \text{and}\quad 0<t<r.
			\end{equation}
			Then we have
			\[ \sum_{i=0}^\infty \rho(\kappa^i r)\apprle \int_0^r \frac{\rho(t)}{t}dt.\]
		\end{enumerate}
	\end{lemma}
	\begin{remark}
		It can be shown that if $f$ satisfies the $\Leb{2}$-$\DMO_x$ in small cylinders and $g$ satisfies the $\Leb{2}$-$\DMO_x$ in small balls, then $\rho_f$ and $\widehat{\rho}_g$ satisfy \eqref{eq:dini-control} (see e.g. \cite[p.495]{Li17}). 
	\end{remark}

	We also recall the following solvability of very weak solutions to parabolic equations whose leading coefficients depend only on $t$. We call such coefficients as \emph{simple coefficients}. This can be proved by a standard duality argument with the solvability result in \cite[Chapter 2]{K07}, see e.g. \cite[Theorem 5]{K09} or  \cite[Lemma 2]{EM17} for the proof.
	\begin{proposition}\label{prop:heat-simple}
		Let $T\in (0,\infty)$ and $q\in (1,\infty)$.  For every $\boldF \in \Leb{q}((0,T)\times\mathbb{R}^d)$, there exists a unique very weak solution $u\in \Leb{q}((0,T)\times\mathbb{R}^d)$ satisfying $u(0,\cdot)=0$ and 
		\begin{equation*}\label{eq:very-weak-sol}
			\partial_t u -a^{ij}(t)D_{ij} u =\Div\Div \boldF\quad \text{in } (0,T)\times\mathbb{R}^d,
		\end{equation*}
		i.e., 
		\begin{equation*}\label{eq:very-weak-sol-weak-form}
			-\int_{(0,T)\times\mathbb{R}^d} u(\partial_t \phi +a^{ij}(t)D_{ij}\phi)\myd{x}dt=\int_{(0,T)\times\mathbb{R}^d} \boldF:\nabla^2 \phi \myd{x}dt
		\end{equation*}
		for all $\phi \in C_0^\infty([0,T)\times\mathbb{R}^d)$. Moreover, we have
		\[ \norm{u}{\Leb{q}((0,T)\times\mathbb{R}^d)}\leq N \norm{\boldF}{\Leb{q}((0,T)\times\mathbb{R}^d)} \]
		for some constant $N=N(d,q,\nu)>0$.
	\end{proposition}

	\subsection{Caccioppoli estimates for Stokes equations}
	In this subsection, we recall Caccioppoli estimates for Stokes equations with simple coefficients.
	\begin{proposition}\label{prop:Caccioppoli-Stokes}
		{Let $0<r<R< 1$.}
		\begin{enumerate}[label=\textnormal{(\roman*)}]
			\item If $u\in \Leb{\infty}\Leb{2}(Q_1)\cap\Leb{2}\tSob{1}(Q_1)$ is a weak solution of 
			\[ \partial_t u-D_i(a^{ij}(t)D_j u)+\nabla \pi=0\quad \text{and}\quad \Div u=g(t)\quad \text{in } Q_1,\]
			for some measurable function $g(t):(-1,0)\rightarrow\mathbb{R}$, then 
			\[  \norm{Du}{\Leb{2}(Q_{r})}\leq \frac{N}{R-r}\norm{u-[u]_{B_R}(t)}{\Leb{2}(Q_R)} \]
			for some constant $N=N(d,\nu)>0$. 
			\item If $u \in \Leb{\infty}\Leb{2}(Q_1)\cap \Leb{2}\tSob{1}(Q_1)$ is a weak solution of 
			\[ \partial_t u - D_i(a^{ij}(t)D_j u)+\nabla \pi =f(t)\quad \text{and}\quad \Div u = g(t)+b(t)\cdot x\quad \text{in } Q_1\]
			for some function $g(t)$ and vector fields $f(t)$, $b(t)$, then
			\[ \norm{D^2u}{\Leb{2}(Q_r)}\leq \frac{N}{(R-r)^2} \norm{u}{\Leb{2}(Q_R)} \]
			for some constant $N=N(d,\nu)>0$. 
		\end{enumerate}
	\end{proposition}
	\begin{proof} 
		(i) See e.g. Dong-Phan \cite[Lemma 2.5]{DP21} and Dong-Kim-Phan \cite[Lemma 3.1]{DKP22} for the proof.
		
		(ii) {By a method of finite difference, one can show that $D^2 u$ satisfies
		\[ \partial_t (D^2 u)-D_i(a^{ij}(t)D_j(D^2 u))+\nabla (D^2 \pi)=0\quad \text{and}\quad \Div (D^2 u)=0\quad \text{in } Q_{7/8},\]
		it follows from (i) that }
		\begin{equation}\label{eq:D3-estimate}
			\norm{D^3 u}{\Leb{2}(Q_r)}\leq \frac{N}{R-r}\norm{D^2 u}{\Leb{2}(Q_R)}
		\end{equation}
		for some constant $N=N(d,\nu)>0$. By an interpolation inequality, \eqref{eq:D3-estimate}, and Young's inequality, for $\varepsilon>0$, we have 
		\begin{equation*}
			\begin{aligned}
				\norm{D^2 u}{\Leb{2}(Q_r)}&\leq \varepsilon \norm{D^3 u}{\Leb{2}(Q_r)}+N\left(\frac{1}{\varepsilon^2}+\frac{1}{r^2}\right) \norm{u}{\Leb{2}(Q_r)}\\
				&\leq \frac{N\varepsilon}{R-r} \norm{D^2 u}{\Leb{2}(Q_R)}+N\left(\frac{1}{\varepsilon^2}+\frac{1}{r^2}\right) \norm{u}{\Leb{2}(Q_R)},
			\end{aligned}
		\end{equation*}
		which implies
		\begin{equation*}
			\norm{D^2 u}{\Leb{2}(Q_r)}\leq \varepsilon  \norm{D^2 u}{\Leb{2}(Q_R)}+\frac{N}{(R-r)^2}\left(1+\frac{1}{\varepsilon^2}\right) \norm{u}{\Leb{2}(Q_R)}
		\end{equation*}
		for some constant $N=N(d,\nu)>0$. Then the desired estimate follows by a standard iteration argument (see e.g. \cite[Lemma 1]{DK11} or \eqref{eq:Du-estimate-preparation}).
	\end{proof}
	
	\section{Estimating the decay rate of the mean oscillation of vorticity}\label{sec:vorticity-rate}
	In this section, we estimate the convergence rate of the mean oscillation of the vorticity and its gradients.  For $X_0=(t_0,x_0)$ and $r>0$, we define 
	\[ \psi(\omega,r,X_0)=\left(\fint_{Q_r(X_0)} |\omega-(\omega)_{Q_r(X_0)}|^2 \myd{x}dt\right)^{1/2},\quad \psi(\omega,r)=\psi(\omega,r,0).\]
	
	The following proposition will play crucial roles in obtaining $C^1_x$ and $C^2_x$-estimates for Stokes equations in divergence form and nondivergence form, respectively.
	\begin{proposition}\label{prop:vorticity-rate}
		Let $r\in (0,1/4)$, $\beta \in (0,1)$, and $X_0 \in Q_{3/4}$. There exists a constant $\kappa \in (0,1/2)$ depending only on $d$, $\nu$, and $\beta$ such that the following hold:
		\begin{enumerate}[label=\textnormal{(\roman*)}]
			\item If $u\in \Leb{\infty}\Leb{2}(Q_1)\cap \Leb{2}\tSob{1}(Q_1)$ is a weak solution of \eqref{eq:div-form} for some $\boldF \in \Leb{2}(Q_1)$, then for any $j\in\mathbb{N}$, we have
			\begin{equation}\label{eq:vorticity-Dini-estimate}
				\psi(\omega,\kappa^j r,X_0)\apprle_{d,\nu,\beta} \kappa^{j\beta}\psi(\omega,r,X_0)+\norm{Du}{\Leb{\infty}(Q_r(X_0))}\Rho_A(\kappa^j r)+\Rho_{\boldF}(\kappa^j r),
			\end{equation}
			where 
			\begin{equation}\label{eq:Rho-defn}
				\Rho_f(r)= \sum_{l=1}^\infty \kappa^{l\beta}(\rho_f(\kappa^{-l}r)[\kappa^{-l}r<1]+\rho_f(1)[\kappa^{-l}r\geq 1]).
			\end{equation}
			Moreover, we have
			\begin{equation}\label{eq:vorticity-summation}
				\sum_{j=0}^\infty \psi(\omega,\kappa^j r, X_0)
				\apprle_{d,\nu,\beta} \psi(\omega,r,X_0)+ \norm{Du}{\Leb{\infty}(Q_r(X_0))} I_{\Rho_A}(r)+I_{\Rho_{\boldF}}(r).
			\end{equation}
			\item If $(u,\pi)\in \tilde{W}^{1,2}_2(Q_1)\times \Sob{0,1}{1}(Q_1)$ is a strong solution of \eqref{eq:nondiv-form} in $Q_1$, then for $j\in\mathbb{N}$,
			\begin{equation}\label{eq:D-vorticity-Dini-estimate}
				\psi(D\omega,\kappa^j r,X_0)
				\apprle_{d,\nu,\beta} \kappa^{j\beta}\psi(D\omega,r,X_0)+\norm{D^2u}{\Leb{\infty}(Q_r(X_0))}\Rho_A(\kappa^j r)+\Rho_{f}(\kappa^j r).
			\end{equation}
			Moreover, we have
			\begin{equation}\label{eq:D-vorticity-summation}
				\sum_{j=0}^\infty \psi(D\omega,\kappa^j r, X_0)
				\apprle_{d,\nu,\beta} \psi(D\omega,r,X_0)+ \norm{D^2u}{\Leb{\infty}(Q_r(X_0))} I_{\Rho_A}(r)+{I_{\Rho_{f}}(r).}
			\end{equation}
		\end{enumerate} 
	\end{proposition}
	\begin{remark}
	By the approximation argument given in the proof of Theorems \ref{thm:A} or \ref{thm:B} later, we may assume that $Du$ or $D^2u$ are bounded.
	\end{remark}

	\begin{proof}
		By translation of the coordinates, we may assume that $X_0=(0,0)$.  For simplicity, we define
		\begin{equation}\label{eq:hat-notation}
		 \hat{a}^{ij}(t)=\fint_{B_r} a^{ij}(t,x)\myd{x}=[a^{ij}]_{0,r}(t)
		\end{equation}
		and similarly, we define
		$\hat{F}^{ij}(t)$.

		(i) Rewrite \eqref{eq:div-form} into
		\[
		\partial_t u -\mathcal{L}_0u+\nabla \pi = \Div \boldF + D_i((a^{ij}-\hat{a}^{ij})D_ju)\quad \text{in } Q_r,
		\]
		where
		\[ \mathcal{L}_0u =D_i(\hat{a}^{ij}(t)D_{j} u).\]
		For $k,l=1,\dots d$, define the vorticity $\omega_{kl}=D_l u^k-D_k u^l$ of the velocity $u$. Then one can see that $\omega_{kl}$ is a very weak solution of
		\[
		\begin{aligned}
			\partial_t\omega_{kl}-\hat{a}^{ij}(t)D_{ij}\omega_{kl}=D_{il}(F^{ki}+(a^{ij}-\hat{a}^{ij})D_ju^k)-D_{ik}(F^{li}+(a^{ij}-\hat{a}^{ij})D_ju^l)
		\end{aligned}
		\]
		in $Q_r$. Since $\hat{a}^{ij}$ depends only on $t$, it follows from Proposition \ref{prop:heat-simple} that there exists a unique very weak solution $\omega_1 \in \Leb{2}((-r^2,0)\times \mathbb{R}^d)$ of
		\[
		\begin{aligned}
			\partial_t (\omega_1)_{kl}-\hat{a}^{ij}D_{ij}(\omega_1)_{kl}&=D_{il}(1_{Q_r}((F^{ki}-\hat{F}^{ki})+(a^{ij}-\hat{a}^{ij})D_j u^k))\\
			&\relphantom{=}-D_{ik}(1_{Q_r}((F^{li}-\hat{F}^{li})+(a^{ij}-\hat{a}^{ij})D_j u^l))
		\end{aligned}
		\]
		satisfying $(\omega_1)_{kl}(-r^2,\cdot)=0$. Moreover, we have
		\begin{equation}\label{eq:omega1-L2}
			\begin{aligned}
				\norm{\omega_1}{\Leb{2}(Q_r)}&\leq N\left(\int_{Q_r} |\boldF-\hat{\boldF}|^2 \myd{x}dt\right)^{1/2}+N\left(\int_{Q_r} |(A-\hat{A})Du|^2 \myd{x}dt\right)^{1/2}
			\end{aligned}
		\end{equation}
		for some constant $N=N(d,\nu)>0$.

		Define $\omega_2=\omega-\omega_1$. Then $\omega_2$ is a very weak solution of 
		\[ \partial_t \omega_2-\hat{a}^{ij}(t)D_{ij}\omega_2=0\quad \text{in }Q_r.\]
		Since $\hat{a}^{ij}$ depends only on $t$, it follows from $C^{1/2,1}_{t,x}$-estimates (see e.g. \cite[Lemma 3]{DK11}) for parabolic equations with simple coefficients and a scaling argument that
		\begin{equation}\label{eq:omega2-Lipschitz}
			[\omega_2]_{C^{1/2,1}_{t,x}(Q_{r/2})}\leq Nr^{-1}\left(|\omega_2-(\omega_2)_{Q_r}|^2\right)^{1/2}_{Q_r}
		\end{equation}
		holds for some constant $N=N(d,\nu)>0$.
		Hence by \eqref{eq:omega1-L2} and \eqref{eq:omega2-Lipschitz}, for $\kappa \in (0,1/2)$, we have
		\begin{align*}
			&(|\omega-(\omega)_{Q_{\kappa r}}|^2)^{1/2}_{Q_{\kappa r}}\leq 2(|\omega_1|^2)^{1/2}_{Q_{\kappa r}}+(|\omega_2-(\omega_2)_{Q_{\kappa r}}|^2)^{1/2}_{Q_{\kappa r}} \\
			&\leq 2(|\omega_1|^2)^{1/2}_{Q_{\kappa r}}+N_0\kappa r[\omega_2]_{C^{1/2,1}(Q_{\kappa r})}\\
			&\leq N_0\kappa(|\omega-(\omega)_{Q_r}|^2)^{1/2}_{Q_r}+N_0(\kappa+\kappa^{-(d+2)/2})(|\omega_1|^2)^{1/2}_{Q_r}\\
			&\leq N_0\kappa (|\omega-(\omega)_{Q_r}|^2)^{1/2}_{Q_r}+N_0(\kappa +\kappa^{-(d+2)/2})(\rho_A(r)\norm{Du}{\Leb{\infty}(Q_r)}+\rho_{\boldF}(r))
		\end{align*}
		for some constant $N_0=N_0(d,\nu)>0$. Then it follows that 
		\begin{equation}\label{eq:psi-rate}
			\psi(\omega,\kappa r)\leq N_0\kappa \psi(\omega,r)+N_0(\kappa +\kappa^{-(d+2)/2})(\rho_A(r)\norm{Du}{\Leb{\infty}(Q_r)}+\rho_{\boldF}(r))
		\end{equation}
		for some constant $N_0=N_0(d,\nu)>0$.
		
		Now given $\beta \in (0,1)$, choose $\kappa \in (0,1)$ so that $N_0<\kappa^{\beta-1}$. Then by \eqref{eq:psi-rate}, $\psi$ satisfies
		\[ \psi(\omega,\kappa r)\leq \kappa^\beta\psi(\omega,r)+N(\rho_A(r)\norm{Du}{\Leb{\infty}(Q_r)}+\rho_{\boldF}(r)) \]
		for some constant $N=N(d,\nu)>0$. Then by iteration, for $j\in \mathbb{N}$, we get \eqref{eq:vorticity-Dini-estimate}. By summing \eqref{eq:vorticity-Dini-estimate} with respect to $j=1$, $2$, ..., and using Lemma \ref{lem:Dini-estimates} (ii), we get  \eqref{eq:vorticity-summation}. This proves (i).
		
		(ii) Similar to (i), we rewrite \eqref{eq:nondiv-form} into
		\[
		\partial_t u -\mathcal{L}_0u+\nabla \pi = f +(a^{ij}-\hat{a}^{ij})D_{ij}u\quad\text{and}\quad \Div u=g,
		\]
		where 
		\begin{equation}\label{eq:L0-notation}
		 \mathcal{L}_0u= \hat{a}^{ij}D_{ij}u .
		\end{equation}
		Since $\hat{a}^{ij}$ depends only on $t$, it follows that there exists a unique weak solution $\omega_1\in \mathcal{H}^1_2((-r^2,0)\times\mathbb{R}^d)$ to
		\[ \partial_t\omega_1-D_i(\hat{a}^{ij}D_{j}\omega_1) =\nabla \times ( (f-\hat{f})1_{Q_r} +(a^{ij}-\hat{a}^{ij})(D_{ij}u) 1_{Q_r})  \]
		in $(-r^2,0)\times \mathbb{R}^d$ satisfying $\omega_1(-r^2,\cdot)=0$. 
		Moreover, we have
		\begin{equation*}
			\begin{aligned}
				\norm{D\omega_1}{\Leb{2}((-r^2,0)\times \mathbb{R}^d)}&\leq N \left(\int_{Q_r}|f-\hat{f}|^2 \myd{x}dt\right)^{1/2}+N\left(\int_{Q_r} |(A-\hat{A})D^2u|^2 \myd{x}dt \right)^{1/2}
			\end{aligned}
		\end{equation*}
		for some constant $N=N(d,\nu)>0$.
		
		Define $\omega_2=\omega-\omega_1$. Then $\omega_2$ satisfies 
		\[ \partial_t\omega_2-\mathcal{L}_0\omega_2=0\quad \text{in } Q_{2r/3}\]
		and it follows from $C^{1/2,1}_{t,x}$-estimates for parabolic equations with simple coefficients that
		\[  [D\omega_2]_{C^{1/2,1}_{t,x}(Q_{r/3})}\apprle r^{-1} (|D\omega_2-(D\omega_2)_{Q_{2r/3}}|^2)_{Q_{2r/3}}^{1/2}.\]
		Then following exactly the same argument as in (i), we get 
		\[ \psi(D\omega,\kappa r)\leq N_0 \kappa \psi(D\omega,r) + N_0(\rho_A(r)\norm{D^2u}{\Leb{\infty}(Q_r)} + \rho_f(r))\]
		for some constant $N_0=N_0(d,\nu)>0$. 
		
		By choosing $\kappa \in (0,1)$ small so that $N_0<\kappa^{\beta-1}$, we get by induction that \eqref{eq:D-vorticity-Dini-estimate} holds. Also, we get \eqref{eq:D-vorticity-summation} by summing \eqref{eq:D-vorticity-Dini-estimate} with respect to $j=1,2,....$ and using Lemma \ref{lem:Dini-estimates} (ii). This completes the proof of Proposition \ref{prop:vorticity-rate}.
	\end{proof}
	
	\begin{remark}\label{rem:rate-of-psi}
		For $0<r<R$, it follows from \eqref{eq:vorticity-Dini-estimate} that
		\begin{equation}\label{eq:rate-vorticity}
			\psi(\omega,r,X_0)\apprle \left(\frac{r}{R}\right)^\beta \psi(\omega,R,X_0)+\norm{Du}{\Leb{\infty}(Q_R(X_0))} \Rho_A(r) +\Rho_{\boldF}(r).
		\end{equation}
		Indeed, choose $\kappa \in (0,1/2)$ in Proposition \ref{prop:vorticity-rate} and then choose $j$ so that $\kappa^{j+1} <r/R\leq \kappa^j$. If $j=0$, then 
		\[ \psi(\omega,r,X_0)\leq \psi(\omega,R,X_0)\apprle \left(\frac{r}{R}\right)^\beta \psi(\omega,R,X_0),\]
		which implies \eqref{eq:rate-vorticity}. If $j\geq 1$, then by \eqref{eq:vorticity-Dini-estimate}, we have
		\begin{align*}
			\psi(\omega,r,X_0)&\apprle \kappa^{j\beta} \psi(\omega,\kappa^{-j} r,X_0)+\norm{Du}{\Leb{\infty}(Q_{\kappa^{-j}r}(X_0))}\Rho_A(r) + \Rho_{\boldF}(r)\\
			&\apprle \left(\frac{r}{R}\right)^{\beta} \psi(\omega,R,X_0)+\norm{Du}{\Leb{\infty}(Q_R(X_0))}\Rho_A(r)+\Rho_{\boldF}(r).
		\end{align*}
		
		Similarly, it follows from \eqref{eq:D-vorticity-Dini-estimate} that
		\begin{equation}\label{eq:rate-D-vorticity}
			\psi(D\omega,r,X_0)\apprle \left(\frac{r}{R}\right)^\beta \psi(D\omega,R,X_0)+\norm{D^2u}{\Leb{\infty}(Q_R(X_0))} \Rho_A(r) +\Rho_{f}(r).
		\end{equation}
	\end{remark}

	\section{Estimating the convergence rate of approximation to velocity}\label{sec:velocity-rate}
	
	In this section, we prove the key proposition which will be used in the proof of the main theorems. For $X_0=(t_0,x_0)$, introduce the following functions
	\begin{equation*}
		\phi_k(r,X_0)=\sup_{t\in (t_0-r^2,t_0)}r^{-k-d/2}\inf_{p\in \mathcal{P}_k}  \norm{u(t,\cdot)-p}{\Leb{2}(B_r(x_0))},
	\end{equation*}
	where $\mathcal{P}_k$ denotes the set of polynomials in $x$ up to order $k$. 
	
	\begin{proposition}\label{prop:velocity-approximation}
		Let $X_0 \in Q_{3/4}$, $r\in (0,1/4)$, and $\beta\in (0,1)$. Then there exists a constant $\kappa\in (0,1/2)$ depending only on $d$, $\nu$, and $\beta$ such that the following hold:
		\begin{enumerate}[label=\textnormal{(\roman*)}]
			\item If $u \in \Leb{\infty}\Leb{2}(Q_1)\cap\Leb{2}\tSob{1}(Q_1)$ is a weak solution of \eqref{eq:div-form} in $Q_1$, then for $j=1,2,\dots$, we have
			\begin{equation}\label{eq:phi-1-control}
				\begin{aligned} 
					\phi_1(\kappa^j r,X_0)&\apprle_{d,\nu,\beta} \kappa^{j\beta}(\phi_1(r,X_0)+\psi(\omega,r,X_0))\\
					&\relphantom{=}+\norm{Du}{\Leb{\infty}(Q_r(X_0))}\Rho_A(\kappa^j r)+\Rho_{\boldF}(\kappa^j r)+\widehat{\Rho}_g(\kappa^j r).
				\end{aligned}
			\end{equation}
			\item If $(u,\pi)\in\tilde{W}^{1,2}_2(Q_1)\times\Sob{0,1}{1}(Q_1)$ is a strong solution of \eqref{eq:nondiv-form} in $Q_1$, then for $j=1,2,\dots$, we have 
			\begin{equation}\label{eq:phi-2-control}
				\begin{aligned} 
					\phi_2(\kappa^j r,X_0)&\apprle_{d,\nu,\beta} \kappa^{j\beta}(\phi_2(r,X_0)+\psi(D\omega,r,X_0))\\
					&\relphantom{=}+\norm{D^2u}{\Leb{\infty}(Q_r(X_0))}\Rho_A(\kappa^j r)+\Rho_{f}(\kappa^j r)+\widehat{\Rho}_{Dg}(\kappa^j r).
				\end{aligned}
			\end{equation}
		\end{enumerate}
	\end{proposition}

	The remaining part of this section consists of three subsections. To obtain \eqref{eq:phi-1-control} and \eqref{eq:phi-2-control}, we derive a systems for $\tilde{u}=u-q$, where $q$ is a suitably chosen polynomial which approximates $u$. Section \ref{subsec:reduction} is devoted to obtaining such systems. After that, we prove Proposition \ref{prop:velocity-approximation} in Sections \ref{subsec:velocity-1} and \ref{subsec:velocity-2}.

	\subsection{Reduction}\label{subsec:reduction}
	We first observe that there exists a radial function $\varphi \in C_0^\infty(B_1)$ satisfying 
	\begin{equation}\label{eq:radial-function-reduction}
		\int_{B_1} \varphi \myd{x}=1\quad \text{and }\quad \int_{B_1} |x|^2 \varphi(x)\myd{x}=0.
	\end{equation}
	Indeed, if $\tilde{\rho}\in C^\infty_0(B_1)$ is a radial function satisfying $\int_{B_1} \tilde{\rho}\myd{x}=1$, then one can check that 
	\[ \varphi(x)=-\frac{1}{3}\tilde{\rho}(x)+\frac{2^{d+2}}{3}\tilde{\rho}(2x) \]
	satisfies \eqref{eq:radial-function-reduction}.

	The following lemma will be used to estimate the velocity field. 
	\begin{lemma}\label{lem:polynomial-comparison}
		For $r>0$, define $\varphi_r(x)=\varphi(x/r)$ and 
		\begin{equation*}\label{eq:mollification}
			u^{(r)}(t,x)=\frac{1}{r^d}\int_{B_r} \varphi_r(y)u(t,x+y)\myd{y}.
		\end{equation*}
		\begin{enumerate}[label=\textnormal{(\roman*)}]
			\item Define $\tilde{u}=u-q_1$, where 
			\begin{equation}\label{eq:q-1-polynomial}
				q_1=q_1(t,x)=u^{(r)}(t,0)+x_k(D_k u^{(r)})(t,0).
			\end{equation}
			Then there exists a constant $N=N(d)>0$ such that 
			\begin{equation*}\label{eq:q-1-polynomial-estimate}
				\norm{\tilde{u}(t,\cdot)}{\Leb{2}(B_r)}\leq N \inf_{p\in \mathcal{P}_1} \norm{u(t,\cdot)-p}{\Leb{2}(B_r)}
			\end{equation*}
			for $r\in (0,1/4)$ and $t\in (-1,0)$.
			\item Define $\tilde{u}=u-q_2$, where 
			\begin{equation}\label{eq:q-2-polynomial}
				q_2=q_2(t,x)=u^{(r)}(t,0)+x_k(D_k u^{(r)})(t,0)+\frac{1}{2}x_jx_k(D_{jk}u^{(r)})(t,0).
			\end{equation}
			Then there exists a constant $N=N(d)>0$ such that 
			\begin{equation*}\label{eq:q-1-polynomial-estimate}
				\norm{\tilde{u}(t,\cdot)}{\Leb{2}(B_r)}\leq N \inf_{p\in \mathcal{P}_2} \norm{u(t,\cdot)-p}{\Leb{2}(B_r)}
			\end{equation*}
			for $r\in (0,1/4)$ and $t\in (-1,0)$.
		\end{enumerate}
	\end{lemma}

	\begin{proof}
		We only prove (ii) since the case (i) is similar. By properties of mollification and H\"older's inequality, we have
		\begin{align*}
			\norm{{q}_2(t,\cdot)}{\Leb{2}(B_r)}\apprle \norm{u(t,\cdot)}{\Leb{2}(B_r)}
		\end{align*}
		and hence it follows that 
		\begin{align*}
			\norm{u(t,\cdot)-{q}_2(t,\cdot)}{\Leb{2}(B_r)}&\apprle \norm{u(t,\cdot)}{\Leb{2}(B_r)}.
		\end{align*}
		On the other hand, note that if $p\in \mathcal{P}_2$ and  $v(t,x)=u(t,x)-p(x)$, then $\tilde{v} =\tilde{u}$. Indeed, if we write
		\[ p(x)=a_0 + b_i x_i +c_{ij} x_i x_j,\quad \text{where }\, c_{ij}=c_{ji},\]
		then by the choice of $\varphi$ in \eqref{eq:radial-function-reduction}, one can easily show that 
		\[ p^{(r)}(t,0)=a_0,\quad (D_k p^{(r)})(t,0)=b_k,\quad \text{and}\quad (D_{ij} p^{(r)})(t,0)=2c_{ij},\]
		which implies that $\tilde{v}=\tilde{u}$. Hence, we get
		\begin{align*}
			\norm{\tilde{u}(t,\cdot)}{\Leb{2}(B_r)}&\apprle \norm{u(t,\cdot)-p}{\Leb{2}(B_r)}.
		\end{align*}
		for any $p\in \mathcal{P}_2$. By taking the infimum over $p\in\mathcal{P}_2$, we get the desired result.
	\end{proof}
	
	Next, we derive a systems for $\tilde{u}=u-q$, where $q$ is a suitably chosen polynomial which approximates $u$. 
	
	We first consider equations in nondivergence form. If $(u,\pi)$ is a strong solution of \eqref{eq:nondiv-form} in $Q_r$, then  $\tilde{u}=u-q_2$, where $q_2$ is defined in \eqref{eq:q-2-polynomial}, satisfies 
	\begin{equation}\label{eq:tilde-u-nondiv}
		\begin{aligned}
			\partial_t\tilde{u}-\mathcal{L}_0\tilde{u}+\nabla \pi&=f+ (a^{ij}-\hat{a}^{ij})D_{ij}u -\partial_t u^{(r)}(t,0)-x_i(D_i\partial_t u^{(r)})(t,0)\\
			&\relphantom{=}-\frac{1}{2}x_i x_j (D_{ij} \partial_t u^{(r)})(t,0)+\hat{a}^{ij}(D_{ij} u^{(r)})(t,0),\\
			\Div\tilde{u}&=g-[g^{(r)}(t,0)+x\cdot Dg^{(r)}(t,0)]
		\end{aligned}
	\end{equation}
	in $Q_r$. Since $\partial_t u^{(r)}(t,0)$ and $\hat{a}^{ij}(D_{ij}u^{(r)})(t,0)$ do not depend on $x$, these terms are absorbed into the pressure. 
	
	To estimate fourth and fifth terms {on the right-hand side of the first equation} in \eqref{eq:tilde-u-nondiv}, note that 
	\[ \partial_t u^{(r)}=\mathcal{L}_0u^{(r)}-\nabla \pi^{(r)}+f^{(r)}+[(a^{ij}-\hat{a}^{ij})D_{ij}u]^{(r)}.\]
	This gives
	\begin{align*}
		&x_k [D_k\partial_t u^{(r)}](t,0)+x_k x_l [D_{kl}\partial_t u^{(r)}](t,0)\\
		&=x_k (D_k\mathcal{L}_0 u^{(r)})(t,0)-x_k D_k \nabla \pi^{(r)}(t,0)-x_k D_k f^{(r)}(t,0)\\
		&\relphantom{=}+x_k D_k[(a^{ij}-\hat{a}^{ij})D_{ij}u]^{(r)}(t,0)+\frac{1}{2}x_kx_l( D_{kl}\mathcal{L}_0 u^{(r)})(t,0)\\
		&\relphantom{=}-\frac{1}{2}x_kx_l D_{kl} \nabla \pi^{(r)}(t,0)-\frac{1}{2}x_kx_l D_{kl} f^{(r)}(t,0)+\frac{1}{2}x_kx_l D_{kl}[(a^{ij}-\hat{a}^{ij})D_{ij}u]^{(r)}(t,0)\\
		&=:\mathrm{I}_1+\mathrm{I}_2+\mathrm{I}_3+\mathrm{I}_4,
	\end{align*}
	where
	\begin{align*}
		\mathrm{I}_1&=x_k (D_k\mathcal{L}_0u^{(r)})(t,0)+\frac{1}{2}x_k x_l (D_{kl}\mathcal{L}_0u^{(r)})(t,0)\\
		\mathrm{I}_2&=-x_k (D_k \nabla \pi^{(r)})(t,0)-\frac{1}{2}x_k x_l (\nabla D_{kl} \pi^{(r)})(t,0),\\
		\mathrm{I}_3&=-x_k (D_k f^{(r)})(t,0)-\frac{1}{2}x_k x_l (D_{kl} f^{(r)})(t,0),\\
		\mathrm{I}_4&=x_k (D_k[(a^{ij}-\hat{a}^{ij})D_{ij}u]^{(r)})(t,0) + \frac{1}{2}x_k x_l (D_{kl}[(a^{ij}-\hat{a}^{ij})D_{ij} u]^{(r)})(t,0).
	\end{align*}
	
	For $\mathrm{I}_2$, we note that 
	\begin{align*}
		D_i\left(\frac{1}{2} x_k x_m (D_{km}\pi^{(r)})(t,0)\right)&=x_k (D_{ki}\pi^{(r)})(t,0),\\
		D_i\left(\frac{1}{3} x_k x_m x_l (D_{kml}\pi^{(r)})(t,0)\right)&=x_k x_l (D_{ikl}\pi^{(r)})(t,0).
	\end{align*}
	Hence, the $\mathrm{I}_2$ term can be absorbed into the pressure. 
	
	Since 
	\begin{align*}
		D_k f^{(r)}(t,0)&=\frac{1}{r^{d+1}}\int_{B_r}(D_k \rho)\left(\frac{y}{r}\right)[f(t,y)-\hat{f}(t)]\myd{y},\\
		D_{kl} f^{(r)}(t,0)&=\frac{1}{r^{d+2}}\int_{B_r}(D_{kl} \rho)\left(\frac{y}{r}\right)[f(t,y)-\hat{f}(t)]\myd{y},
	\end{align*}
	it follows from H\"older's inequality that 
	\begin{equation}\label{eq:I3-estimate-nondiv}
		\norm{\mathrm{I}_3}{\Leb{2}(Q_r)}\apprle r^{d/2+1} \rho_f(r).
	\end{equation}
	Similarly, we have
	\begin{equation}\label{eq:I4-estimate-nondiv}
		\norm{\mathrm{I}_4}{\Leb{2}(Q_r)}\apprle r^{d/2+1} \rho_A(r) \norm{D^2u}{\Leb{\infty}(Q_r)}.
	\end{equation}
	
	Finally, for $\mathrm{I}_1$, note first that the $j$th component of the first term in $\mathrm{I}_1$ can be rewritten as 
	\begin{equation}\label{eq:vorticity-nondiv-expression}
		x_k (D_k \mathcal{L}_0u_j^{(r)})(t,0)=\frac{1}{2}x_k (\mathcal{L}_0 \omega_{jk}^{(r)})(t,0)+\frac{1}{2}D_j(x_k x_l (\mathcal{L}_0 D_l u_k^{(r)})(t,0)).
	\end{equation}
	To estimate the second term in $\mathrm{I}_1$, we note that 
	\begin{equation}\label{eq:second-vorticity-1}
		\frac{1}{2}x_k x_l (D_{kl} \mathcal{L}_0 u_j^{(r)})(t,0)=\frac{1}{2}x_k x_l (D_k \mathcal{L}_0 \omega_{jl}^{(r)})(t,0)+\frac{1}{2}x_k x_l (D_{kj} \mathcal{L}_0 u_l^{(r)})(t,0).
	\end{equation}

	If we define
	\[ \Phi(t,x)= x_k x_l x_m (D_{kl} \mathcal{L}_0 u_m^{(r)})(t,0),\]
	then one can easily see that 
	\begin{equation}\label{eq:second-vorticity-2}
		D_j \Phi(t,x)=2x_l x_k (D_{jk} \mathcal{L}_0 u_l^{(r)})(t,0)+x_k x_l(D_{kl}\mathcal{L}_0u_j^{(r)})(t,0).
	\end{equation}
	By \eqref{eq:second-vorticity-1} and \eqref{eq:second-vorticity-2}, we have 
	\begin{equation}\label{eq:vorticity-nondiv-expression-2}
		\frac{1}{2}x_k x_l(D_{kl} \mathcal{L}_0u_j^{(r)})(t,0)=\frac{1}{3} x_k x_l (D_k\mathcal{L}_0\omega_{jl}^{(r)})(t,0)+\frac{1}{6}D_j \Phi(t,x).
	\end{equation}
	Hence by \eqref{eq:vorticity-nondiv-expression} and \eqref{eq:vorticity-nondiv-expression-2}, for $j=1,\dots,d$, $\mathrm{I}_1^j$ can be rewritten as 
	\[ 
	\mathrm{I}_1^j=\frac{1}{2}x_k (\mathcal{L}_0\omega_{kj}^{(r)})(t,0)+\frac{1}{3}x_k x_l(D_k\mathcal{L}_0\omega_{lj}^{(r)})(t,0)+D_j\Psi
	\]
	for some $\Psi$. We again absorb $\Psi$ into the pressure and estimate \begin{equation}\label{eq:I1-estimate-nondiv}
		\norm{x_k(\mathcal{L}_0\omega^{(r)})(t,0)}{\Leb{2}(Q_r)}+\norm{x_k x_l(D_k\mathcal{L}_0\omega^{(r)})(t,0)}{\Leb{2}(Q_r)}\apprle r^{d/2+1}\psi(D\omega,r).
	\end{equation}
	Hence for $j=1,\dots d$, equation \eqref{eq:tilde-u-nondiv} can be rewritten as 
	\begin{equation*}\label{eq:tilde-u-nondiv-final}
		\left\{\begin{aligned}
			\partial_t \tilde{u}^j-\mathcal{L}_0\tilde{u}^j+D_j {\tilde{\pi}}&=f^j+(a^{kl}-\hat{a}^{kl})D_{kl}u^j\\
			&\relphantom{=}-\frac{1}{2}x_k(D_k\mathcal{L}_0\omega_{jk}^{(r)})(t,0)-\frac{1}{3}x_k x_l(D_k\mathcal{L}_0\omega_{jl}^{(r)})(t,0)\\
			&\relphantom{=}+\mathrm{I}_3^j +\mathrm{I}_4^j,\\
			\Div \tilde{u}&=g-(g^{(r)}(t,0)+x \cdot (Dg^{(r)})(t,0))
		\end{aligned}
		\right.
	\end{equation*}
	in $Q_r$.

	Similarly, for the equation in divergence form, note that $\tilde{u}=u-q_1$, where $q_1$ is defined in \eqref{eq:q-1-polynomial}, satisfies
	\begin{equation}\label{eq:tilde-u-div}
		\left\{\begin{aligned}
			\partial_t\tilde{u}-\mathcal{L}_0\tilde{u}+\nabla \pi &=\Div \boldF + D_i ((a^{ij}-\hat{a}^{ij})D_j u) -\partial_t u^{(r)}(t,0)\\
			&\relphantom{=}-x\cdot (\nabla \partial_t u^{(r)})(t,0),\\
			\Div\tilde{u}&=g-g^{(r)}(t,0)
		\end{aligned}
		\right.
	\end{equation}
{in $Q_r$.}
	Note that $\partial_t u^{(r)}(t,0)$ is independent of $x$ and can be absorbed into the pressure. 
	
	To estimate fourth term, we note that 
	\[ \partial_t u^{(r)}=\mathcal{L}_0 u^{(r)}-\nabla \pi^{(r)} + D_i ((a^{ij}-\hat{a}^{ij})D_j u)^{(r)} + \Div\boldF^{(r)}.\]
	This gives
	\begin{align*}
		&x\cdot (\nabla \partial_t u^{(r)})(t,0)\\
		&=x\cdot (\nabla \mathcal{L}_0 u^{(r)})(t,0)-x\cdot (D^2 \pi^{(r)})(t,0) + x\cdot (\nabla D_i ((a^{ij}-\hat{a}^{ij})D_j u)^{(r)})(t,0)\\
		&\relphantom{=}+x\cdot (\nabla \Div(\boldF^{(r)}))(t,0):=\mathrm{I}_1 +\mathrm{I}_2 + \mathrm{I}_3 +\mathrm{I}_4.
	\end{align*}
	
	Following a similar argument as in the case of nondivergence form, $\mathrm{I}_2$ can be absorbed into the pressure and 
	\begin{equation}\label{eq:I3-I4-estimate}
		\norm{\mathrm{I}_3}{\Leb{2}(Q_r)}+\norm{\mathrm{I}_4}{\Leb{2}(Q_r)}\leq Nr^{d/2}\left(\rho_A(r)\norm{Du}{\Leb{\infty}(Q_r)}+\rho_{\boldF}(r)\right)
	\end{equation}
	holds for some constant $N=N(d)>0$. Finally, note that 
	\[ \mathrm{I}_1^k = D_k \Phi +\frac{1}{2} x_i(\mathcal{L}_0\omega_{ki}^{(r)})(t,0),\quad j=1,\dots,d,\]
	where 
	\[ \Phi(t,x)=\frac{1}{2}x_i x_j (D_i \mathcal{L}_0 u_j^{(r)})(t,0).\]
	We again absorb $\Phi$ into the pressure and estimate 
	\begin{equation}\label{eq:I1-estimate}
		\norm{x_i(\mathcal{L}_0\omega_{ki}^{(r)})(t,0)}{\Leb{2}(Q_r)}\apprle r^{d/2}\psi(\omega,r).
	\end{equation}
	Hence the equation \eqref{eq:tilde-u-div} can be rewritten as 
	\begin{equation}\label{eq:tilde-u-div-final}
		\left\{\begin{aligned}
			\partial_t\tilde{u}-\mathcal{L}_0\tilde{u}+\nabla {\tilde{\pi}} &=\Div (\boldF-\hat{\boldF})+D_i((a^{ij}-\hat{a}^{ij})D_{j}u)\\
			&\relphantom{=} -\frac{1}{2}x(\mathcal{L}_0\omega^{(r)})(t,0)+\mathrm{I}_3+\mathrm{I}_4,\\
			\Div \tilde{u}&=g-g^{(r)}(t,0)
		\end{aligned}
		\right.
	\end{equation}
	in $Q_r$.
	
	\subsection{Proof of Proposition \ref{prop:velocity-approximation} (i)}\label{subsec:velocity-1}
	In this subsection, we write $\psi(r,X_0)=\psi(\omega,r,X_0)$. First, assume that $X_0=(0,0)$. By Lemma \ref{lem:approximate-system-existence}, there exists $v_1\in \Leb{\infty}\Leb{2}((-r^2,0)\times\mathbb{R}^d)\cap \Leb{2}\tSob{1}((-r^2,0)\times\mathbb{R}^d)$ satisfying 
	\begin{equation}\label{eq:tilde-u-div-final-whole}
		\left\{\begin{alignedat}{2}
			\partial_t v_1-\mathcal{L}_0v_1+\nabla \pi_1 &=\Div ((\boldF-\hat{\boldF}))+D_i((a^{ij}-\hat{a}^{ij})(D_{j}u))\\
			&\relphantom{=} \left( -\frac{1}{2}x(\mathcal{L}_0\omega^{(r)})(t,0)+\mathrm{I}_3+\mathrm{I}_4\right)&&\quad \text{in } {Q_r},\\
			\Div v_1&=g-[g]_{B_r}(t)&&\quad \text{in } Q_r,\\
			v_1(-r^2,\cdot)&=0&&\quad \text{on } \mathbb{R}^d
		\end{alignedat}
		\right.
	\end{equation}
	and
	\begin{equation}\label{eq:v-estimate-div}
		\begin{aligned}
			&\norm{v_1}{\Leb{\infty}\Leb{2}((-r^2,0)\times\mathbb{R}^d)}
			\apprle \norm{\boldF-\hat{\boldF}}{\Leb{2}(Q_r)}+\norm{(a^{ij}-\hat{a}^{ij})Du}{\Leb{2}(Q_r)}\\
			&\relphantom{=}+r\left(\norm{x(\mathcal{L}_0\omega^{(r)})(t,0)}{\Leb{2}(Q_r)}+\norm{\mathrm{I}_3}{\Leb{2}(Q_r)}+\norm{\mathrm{I}_4}{\Leb{2}(Q_r)}\right)+r\norm{g-[g]_{B_r}(t)}{\Leb{\infty}\Leb{2}(Q_r)}.
		\end{aligned}
	\end{equation}
	
	By H\"older's inequality, we have
	\begin{equation}\label{eq:g-estimate}
		\begin{aligned}
			\norm{g-[g]_{B_r}(t)}{\Leb{\infty}\Leb{2}(Q_r)}&\apprle \sup_{t\in (-r^2,0)} \left(\int_{B_r}\fint_{B_r} |g(t,x)-g(t,y)|^2 dydx \right)^{1/2}\\
			&\apprle r^{d/2} \widehat{\rho}_g(r).
		\end{aligned}
	\end{equation}
	
	Hence by  \eqref{eq:I3-I4-estimate}, \eqref{eq:I1-estimate}, \eqref{eq:v-estimate-div}, and \eqref{eq:g-estimate}, we have
	\begin{equation}\label{eq:v1-control-estimates}
		\norm{v_1}{\Leb{\infty}\Leb{2}((-r^2,0)\times\mathbb{R}^d)}\apprle r^{d/2+1}\left(\rho_{\boldF}(r)+\rho_A(r)\norm{Du}{\Leb{\infty}(Q_r)}+\psi(r)+\widehat{\rho}_g(r)\right).
	\end{equation}

	Define $v_2=\tilde{u}-v_1=u-q_1-v_1$. Then by \eqref{eq:tilde-u-div-final} and \eqref{eq:tilde-u-div-final-whole}, $v_2$ satisfies
	\begin{equation*}\label{eq:w-equation}
		\partial_t v_2-\mathcal{L}_0v_2+\nabla \pi_2 = 0\quad \text{and}\quad \Div v_2=[g]_{B_r}(t)-g^{(r)}(t,0)\quad \text{in } Q_{r}.
	\end{equation*}
	Then it follows from Proposition \ref{prop:Caccioppoli-Stokes} that 
	\begin{equation}\label{eq:v2-estimate}
		\norm{Dv_2}{\Leb{2}(Q_{2r/3})}\leq N \norm{v_2}{\Leb{\infty}\Leb{2}(Q_r)}
	\end{equation}
	for some constant $N=N(d,\nu)>0$.  If we define  $\omega_2=\nabla\times v_2$, then {by applying a method of finite difference, one can show that $\omega_2$ satisfies 
	\[ \partial_t \omega_2-\mathcal{L}_0 \omega_2=0\quad \text{in } Q_{r/2}.\]
	Since $\hat{a}^{ij}$ depends only on $t$,} it follows from $C^{\alpha/2,\alpha}_{t,x}$-estimates for such parabolic equations (see e.g. \cite[Lemma 3]{DK11}) and a scaling argument that
	\begin{equation}\label{eq:omega2-estimate}
		[\omega_2]_{C^{\alpha/2,\alpha}_{t,x}(Q_{2r/5})}\leq Nr^{-\alpha}(|\omega_2-(\omega_2)_{Q_{r/2}}|^2)^{1/2}_{Q_{r/2}} 
	\end{equation}
	for some constant $N=N(d,\nu)>0$. Hence by \eqref{eq:v2-estimate} and \eqref{eq:omega2-estimate}, we have
	\begin{equation}\label{eq:omega-2-Holder}
		\begin{aligned}
			[\omega_2]_{C^{\alpha/2,\alpha}_{t,x}(Q_{2r/5})}&\leq Nr^{-\alpha-(d+2)/2}\norm{Dv_2}{\Leb{2}(Q_{r/2})}\\
			&\leq Nr^{-\alpha-1-d/2}\norm{v_2}{\Leb{\infty}\Leb{2}(Q_{2r/3})}
		\end{aligned}
	\end{equation}
	for some constant $N=N(d,\nu,\alpha)>0$. 
	
	On the other hand, we have
	\[ \Delta v^j_2 = D_i\omega^{ji}_2\quad \text{in } Q_r,\quad j=1,\dots,d.\]
	Fix $\beta \in (0,1)$ and choose $\alpha \in (\beta,1)$. Then by the Schauder estimates for the Poisson equation, we have 
	\begin{equation}\label{eq:D2-estimate}
		[Dv_2(t,\cdot)]_{C^{0,\alpha}(B_{r/4})}\leq N\left(r^{-1-\alpha-d/2}\norm{v_2(t,\cdot)}{\Leb{2}(B_{r/3})}+[\omega_2(t,\cdot)]_{C^{0,\alpha}(B_{r/3})} \right)
	\end{equation}
	for some constant $N=N(d,\alpha)>0$. 
	
	Let $\hat{q}_1=\hat{q}_1(t,x)$ be the first-order Taylor expansion of $v_2$ with respect to $x$ at $(t,0)$. Then by \eqref{eq:omega-2-Holder} and \eqref{eq:D2-estimate}, for $\kappa \in (0,1)$, we have
	\begin{equation}\label{eq:v2-Taylor}
		\norm{v_2-\hat{q}_1}{\Leb{\infty}\Leb{2}(Q_{\kappa r})}\leq N \kappa^{1+\alpha+d/2}\norm{v_2}{\Leb{\infty}\Leb{2}(Q_r)}
	\end{equation}
	for some constant $N=N(d,\nu,\alpha)>0$.
	
	Define $q=q_1+\hat{q}_1$, where $q_1$ is defined in \eqref{eq:q-1-polynomial}. Then by \eqref{eq:v1-control-estimates} and \eqref{eq:v2-Taylor}, we have
	\begin{equation}\label{eq:preparation-for-iteration}
		\begin{aligned}
			&\norm{u-q}{\Leb{\infty}\Leb{2}(Q_{\kappa r})}\leq \norm{v_1}{\Leb{\infty}\Leb{2}(Q_{\kappa r})}+\norm{v_2-\hat{q}_1}{\Leb{\infty}\Leb{2}(Q_{\kappa r})}\\
			&\leq Nr^{d/2+1}\left(\rho_A(r)\norm{Du}{\Leb{\infty}(Q_r)}+\rho_{\boldF}(r)+\widehat{\rho}_g(r)+\psi(r)\right)\\
			&\relphantom{=}+N\kappa^{1+\alpha+d/2}\left(\norm{\tilde{u}}{\Leb{\infty}\Leb{2}(Q_r)}+\norm{v_1}{\Leb{\infty}\Leb{2}(Q_r)}\right)\\
			&\leq Nr^{d/2+1}\left(\rho_A(r)\norm{Du}{\Leb{\infty}(Q_r)}+\rho_{\boldF}(r)+\widehat{\rho}_g(r)+\psi(r)\right)\\
			&\relphantom{=}+N\kappa^{1+\alpha+d/2}\norm{\tilde{u}}{\Leb{\infty}\Leb{2}(Q_{r})}
		\end{aligned}
	\end{equation}
	for some constant $N=N(d,\nu,\alpha)>0$. Then by  \eqref{eq:preparation-for-iteration} and Lemma \ref{lem:polynomial-comparison} (i), we have
	\begin{equation}\label{eq:prepare-for-iteration-A}
		\begin{aligned}
			&\sup_{t\in (-(\kappa r)^2,0)} \inf_{p\in\mathcal{P}_1} (\kappa r)^{-1-d/2}\norm{u(t,\cdot)-p}{\Leb{2}(B_{\kappa r})}\\
			&\leq N_1 \kappa^\alpha \sup_{t\in (-r^2,0)} \inf_{p\in\mathcal{P}_1} r^{-1-d/2}\norm{u(t,\cdot)-p}{\Leb{2}(B_r)}\\
			&\relphantom{=}+N\kappa^{-1-d/2} \left(\rho_A(r)\norm{Du}{\Leb{\infty}(Q_r)}+\rho_{\boldF}(r)+\widehat{\rho}_g(r)+\psi(r)\right)
		\end{aligned}
	\end{equation}
	for some constants $N_1$ and $N$ depending only on $d$, $\nu$, and $\alpha$.
	
	Choose $\beta_1\in(\beta,\alpha)$ and a sufficiently small $\kappa\in (0,1)$ so that $N_1\leq \kappa^{\beta_1-\alpha}$. By definition of $\phi_1$ we have 
	\begin{equation}\label{eq:phi-estimate}
		\phi_1(r) \apprle \frac{1}{r^{1+d/2}} \norm{u}{\Leb{\infty}\Leb{2}(Q_r)}.
	\end{equation}
	Then \eqref{eq:prepare-for-iteration-A} is rewritten as 
	\[
	\phi_1(\kappa r)\leq \kappa^{\beta_1} \phi_1(r)+N\kappa^{-1-d/2} \left(\rho_A(r)\norm{Du}{\Leb{\infty}(Q_r)} + \rho_{\boldF}(r)+\widehat{\rho}_g(r)+\psi(r)\right).
	\]
	By iteration, we have 
	\begin{equation}\label{eq:phi-1-iteration}
		\begin{aligned}
			\phi_1(\kappa^j r)
			&\leq \kappa^{j\beta_1} \phi_1(r)+N\kappa^{-1-d/2}\norm{Du}{\Leb{\infty}(Q_r)}\left(\sum_{l=1}^j \kappa^{\beta_1(l-1)} \rho_A(\kappa^{j-l} r)\right)\\
			&\relphantom{=}+N\kappa^{-1-d/2}\left(\sum_{l=1}^j \kappa^{\beta_1(l-1)} \left(\rho_{\boldF}(\kappa^{j-l}r)+\widehat{\rho}_g(\kappa^{j-l}r)+\psi(\kappa^{j-l}r) \right)\right).
		\end{aligned}
	\end{equation}
	
	For $j=1$, we have
	\begin{align*}
		\phi_1(\kappa r)&\leq \kappa^{\beta_1}\phi_1(r)+N\kappa^{-1-d/2}(\norm{Du}{\Leb{\infty}(Q_r)}\rho_A(r)+\rho_{\boldF}(r)+\widehat{\rho}_g(r)+\psi(r)).
	\end{align*}
	
	For $j\geq 2$, by \eqref{eq:vorticity-Dini-estimate}, we note that 
	\begin{align*}
		\psi(\kappa^{j-l} r)&\leq \kappa^{(j-l)\beta}\psi(r)\\
		&\relphantom{=}+N\sum_{m=1}^{j-l} \kappa^{(m-1)\beta}\left(\rho_A(\kappa^{j-l-m}r)\norm{Du}{\Leb{\infty}(Q_r)}+\rho_{\boldF}(\kappa^{j-l-m}r)\right).
	\end{align*}
	Then it follows that 
	\begin{align}\label{eq:psi-1-iteration} 
		&\sum_{l=1}^j \kappa^{(l-1)\beta_1} \psi(\kappa^{j-l}r)\nonumber\\
		&\leq \sum_{l=1}^j \kappa^{(l-1)\beta_1}\kappa^{(j-l)\beta}\psi(r)+N\norm{Du}{\Leb{\infty}(Q_r)}\sum_{l=1}^j \sum_{m=1}^{j-l} \kappa^{(l-1)\beta_1}\kappa^{(m-1)\beta}\rho_A(\kappa^{j-l-m}r)\\
		&\relphantom{=}+N\sum_{l=1}^j \sum_{m=1}^{j-l} \kappa^{(l-1)\beta_1}\kappa^{(m-1)\beta}\rho_{\boldF}(\kappa^{j-l-m}r)\nonumber\\
		&\leq \frac{1}{\kappa^{\beta}-\kappa^{\beta_1}}\left[\kappa^{j\beta} \psi(r)+N\norm{Du}{\Leb{\infty}(Q_r)}\sum_{l=1}^{j} \kappa^{(l-1)\beta} (\rho_A(\kappa^{j-l}r)+ \rho_{\boldF}(\kappa^{j-l}r))\right]. \nonumber
	\end{align}
	Here we used the fact that $\kappa \in (0,1)$, $\beta <\beta_1$, 
	\[  \sum_{l=1}^j \kappa^{(l-1)\beta_1} \kappa^{(j-l)\beta}\leq \kappa^{j\beta} \sum_{l=1}^\infty \kappa^{l(\beta_1-\beta)-\beta_1}=\frac{\kappa^{j\beta}}{\kappa^{\beta}-\kappa^{\beta_1}}\]
	and 
	\begin{align*}
		\sum_{l=1}^j\sum_{m=1}^{j-l} \kappa^{(l-1)\beta_1}\kappa^{(m-1)\beta}\rho_A(\kappa^{j-(l+m)})&\leq\sum_{s=1}^j \sum_{m+l=s} \kappa^{(l-1)\beta_1}\kappa^{(m-1)\beta}\rho_A(\kappa^{j-s}r)\\
		&\leq \frac{1}{\kappa^{\beta}-\kappa^{\beta_1}} \sum_{s=1}^j \kappa^{(s-1)\beta} \rho_A(\kappa^{j-s}r).
	\end{align*}
	
	Hence by \eqref{eq:phi-1-iteration} and \eqref{eq:psi-1-iteration}, we get
	\begin{equation*}\label{eq:kappa-j-r-summation}
		\begin{aligned}
			\phi_1(\kappa^j r)&\apprle_{d,\nu,\beta} \kappa^{j\beta}(\phi_1(r) +\psi(r))+\norm{Du}{\Leb{\infty}(Q_r)}\left(\sum_{l=1}^j \kappa^{\beta(l-1)} \rho_A(\kappa^{j-l}r)\right)\\
			&\relphantom{=}+N\sum_{l=1}^{j} \kappa^{\beta(l-1)}(\rho_{\boldF}(\kappa^{j-l}r)+\widehat{\rho}_g(\kappa^{j-l}r))\\
			&\apprle_{d,\nu,\beta} \kappa^{j\beta}( \phi_1(r)+\psi(r)) +\norm{Du}{\Leb{\infty}(Q_r)}\Rho_A(\kappa^j r)+\Rho_{\boldF}(\kappa^j r)+\widehat{\Rho}_g(\kappa^j r),
		\end{aligned}
	\end{equation*}
	where $\Rho_f$ is given in \eqref{eq:Rho-defn} and 
	\begin{equation}\label{eq:Rho-hat-g-defn}
		\widehat{\Rho}_g(r)=\sum_{l=0}^\infty \kappa^{\beta l}(\widehat{\rho}_{g}(\kappa^{-l}r)[\kappa^{-l}r<1]+\widehat{\rho}_{g}(1)[\kappa^{-l}r\geq 1]).
	\end{equation}
	
	Hence by translation of the coordinates, for any $X_0=(t_0,x_0)\in Q_{3/4}$, we have
	\begin{equation*}\label{eq:kappa-iteration}
		\begin{aligned}
			\phi_1(\kappa^j r,X_0)&\apprle_{d,\nu,\beta} \kappa^{j\beta} (\phi_1(r,X_0)+\psi(r,X_0)) \\
			&\relphantom{=}+\norm{Du}{\Leb{\infty}(Q_r(X_0))}\Rho_A(\kappa^j r)+\Rho_{\boldF}(\kappa^j r)+\widehat{\Rho}_g(\kappa^j r)
		\end{aligned}
	\end{equation*}
	and 
	\begin{equation*}\label{eq:kappa-summation}
		\begin{aligned}
			\sum_{j=0}^\infty\phi_1(\kappa^j r,X_0)&\apprle_{d,\nu,\beta}  \phi_1(r,X_0)+\psi(r,X_0) +\norm{Du}{\Leb{\infty}(Q_r(X_0))}I_{\Rho_A}(r)\\
			&\relphantom{=}+I_{\Rho_{\boldF}}(r)+I_{\widehat{\Rho}_{g}}(r).
		\end{aligned}
	\end{equation*}
	
	This completes the proof of Proposition \ref{prop:velocity-approximation} (i). \hfill\qed
	
	\begin{remark}
		Following a similar argument as in Remark \ref{rem:rate-of-psi}, for $0<s\leq r$, we have
		\begin{equation}\label{eq:phi-psi-iteration-final}
			\begin{aligned}
				\phi_1(s,X_0)&\apprle \left(\frac{s}{r}\right)^\beta \left(\phi_1(r,X_0)+ \psi(r,X_0)\right)\\
				&\relphantom{=}+\norm{Du}{\Leb{\infty}(Q_r(X_0))}\Rho_A(s)+\Rho_{\boldF}(s)+\widehat{\Rho}_g(s).
			\end{aligned}
		\end{equation}
	\end{remark}
	
	\subsection{Proof of Proposition \ref{prop:velocity-approximation} (ii)}\label{subsec:velocity-2} {Recall the decomposition given in Section \ref{subsec:reduction} and the notation $\hat{f}$ and $\mathcal{L}_0v=\hat{a}^{ij}(t)D_{ij}v$ given in \eqref{eq:hat-notation} and \eqref{eq:L0-notation}, respectively.}
	For simplicity, we write 
	\[ \psi(r,X_0)=\psi(D\omega,r,X_0),\quad X_0=(t_0,x_0). \]
	
	We first assume that $X_0=(0,0)$. For $j=1,\dots,d$, define 
	\begin{align*}
		h^j&=\left(f^j-\hat{f}^j+(a^{kl}-\hat{a}^{kl})D_{kl}u^j\right)\\
		&\relphantom{=}-\left(\frac{1}{2}x_k(D_k\mathcal{L}_0\omega_{kj}^{(r)})(t,0)+\frac{1}{3}x_k x_l(D_k\mathcal{L}_0\omega_{lj}^{(r)})(t,0)\right)+\mathrm{I}_3^j +\mathrm{I}_4^j.
	\end{align*}
	By Lemma \ref{lem:approximate-system-existence}, there exists $v_1\in \Leb{\infty}\Leb{2}((-r^2,0)\times\mathbb{R}^d)\cap \Leb{2}\tSob{1}((-r^2,0)\times\mathbb{R}^d)$ satisfying $v_1(-r^2,\cdot)=0$, 
	\begin{equation*}\label{eq:tilde-u-nondiv-final-decom}
		\left\{\begin{alignedat}{2}
			\partial_t v^j_1-\mathcal{L}_0v^j_1+D_j {\tilde{\pi}}&=h^j&&\quad \text{in } {Q_r},\\
			\Div v_1&=g-([g]_{B_r}(t)+x\cdot [Dg]_{B_r}(t))&&\quad \text{in } Q_r,
		\end{alignedat}
		\right.
	\end{equation*}
	and 
	\begin{equation}\label{eq:v1-estimate-nondiv-a}
		\begin{aligned}
			\norm{v_1}{\Leb{\infty}\Leb{2}((-r^2,0)\times\mathbb{R}^d)}
			&\apprle r^{d/2+2}\left(\psi(r)+\rho_f(r)+\rho_A(r)\norm{D^2u}{\Leb{\infty}(Q_r)}\right)\\
			&\relphantom{=}+r\norm{g-[g]_{B_r}(t)-x\cdot[Dg]_{B_r}(t)}{\Leb{\infty}\Leb{2}(Q_r)}.
		\end{aligned}
	\end{equation}

	Here we used \eqref{eq:I3-estimate-nondiv}, \eqref{eq:I4-estimate-nondiv}, and \eqref{eq:I1-estimate-nondiv}. By the Poincar\'e inequality, we have 
	\begin{equation}\label{eq:Dg-estimate-oscillation}
		\begin{aligned}
			&r\norm{g-[g]_{B_r}(t)-x\cdot[Dg]_{B_r}(t)}{\Leb{\infty}\Leb{2}(Q_r)}\\
			&\apprle r^2\norm{Dg-[Dg]_{B_r}(t)}{\Leb{\infty}\Leb{2}(Q_r)}\apprle r^{d/2+2}\hat{\rho}_{Dg}(r).
		\end{aligned}
	\end{equation}
	Hence by \eqref{eq:v1-estimate-nondiv-a} and \eqref{eq:Dg-estimate-oscillation}, we get
	\begin{equation}\label{eq:v1-estimate-nondiv}
		\norm{v_1}{\Leb{\infty}\Leb{2}((-r^2,0)\times\mathbb{R}^d)}\apprle r^{d/2+2}\left(\psi(r)+\rho_f(r)+\rho_A(r)\norm{D^2u}{\Leb{\infty}(Q_r)}+\widehat{\rho}_{Dg}(r)\right).
	\end{equation}

	Define $v_2:=\tilde{u}-v_1=u-q_2-v_1$, where $q_2$ is defined in \eqref{eq:q-2-polynomial}. Then $v_2$ satisfies
	\[ \partial_t v_2 -\mathcal{L}_0 v_2 +\nabla \pi_2 = \hat{f}(t)\quad \text{in } Q_{r} \]
	and 
	\begin{equation}\label{eq:v2-div-nonfree}
		\Div v_2 =[g]_{B_r}(t)+x\cdot [Dg]_{B_r}(t)-\left(g^{(r)}(t,0)+x\cdot Dg^{(r)}(t,0)\right)\quad \text{in } Q_{r}.
	\end{equation}
	
	Note that $\omega_2:=\nabla \times v_2$ satisfies
	\begin{equation}\label{eq:omega2-equation}
		\partial_t \omega_2-\mathcal{L}_0\omega_2=0\quad \text{in } Q_{r/2}.
	\end{equation}
	
	Since $\hat{a}^{ij}$ depends only on $t$ and $D\omega_2$ is a solution of \eqref{eq:omega2-equation}, it follows from $C^{\alpha/2,\alpha}_{t,x}$-estimates of \eqref{eq:omega2-equation} and a scaling argument that
	\begin{equation}\label{eq:Domega-2-Holder}
		[D\omega_2]_{C^{\alpha/2,\alpha}_{t,x}(Q_{r/3})}\leq N r^{-\alpha} (|D\omega_2-(D\omega_2)_{Q_{2r/3}}|^2)_{Q_{2r/3}}^{1/2}.
	\end{equation}
	
	Since $\hat{a}^{ij}$ depends only on $t$ and $v_2$ satisfies \eqref{eq:v2-div-nonfree}, it follows from Proposition \ref{prop:Caccioppoli-Stokes} (ii) that 
	\begin{equation}\label{eq:Hessian-u2}
		\norm{D^2v_2}{\Leb{2}(Q_{2r/3})}\leq Nr^{-1} \norm{v_2}{\Leb{\infty}\Leb{2}(Q_r)}
	\end{equation}
	for some constant $N=N(d,\nu)>0$.
	
	Fix $\beta \in (0,1)$ and choose $\alpha \in (\beta,1)$. Then note that 
	\[ \Delta v_2^j = D_i (\omega_2^{ji})+[D_jg]_{B_r}(t)-(D_jg)^{(r)}(t,0),\quad j=1,\dots,d.\]
	Then it follows from an interior Schauder estimates for the Poisson equation  that 
	\begin{equation}\label{eq:C2alpha-Laplacian}
		[D^2 v_2(t,\cdot)]_{C^{0,\alpha}(B_{r/4})}\leq N\left(r^{-2-\alpha-d/2}\norm{v_2(t,\cdot)}{\Leb{2}(B_r)}+[D\omega_2(t,\cdot)]_{C^{0,\alpha}(B_{r/3})} \right)
	\end{equation}
	for some constant $N=N(d,\alpha)>0$.

	By \eqref{eq:Domega-2-Holder} and \eqref{eq:Hessian-u2}, we have
	\begin{equation}\label{eq:Domega-estimates}
		\begin{aligned}
			[D\omega_2]_{C^{\alpha/2,\alpha}_{t,x}(Q_{r/3})}&\leq Nr^{-\alpha}(|D\omega_2-(D\omega_2)_{Q_{2r/3}}|^2)^{1/2}_{Q_{2r/3}}\\
			&\leq N r^{-\alpha-1-d/2} \norm{D^2 v_2}{\Leb{2}(Q_{2r/3})}\\
			&\leq N r^{-\alpha-2-d/2} \norm{v_2}{\Leb{\infty}\Leb{2}(Q_r)}
		\end{aligned}
	\end{equation}
	for some constant $N=N(d,\nu,\alpha)>0$.
	
	Let $\hat{q}_2=\hat{q}_2(t,x)$ be the second-order Taylor expansion of $v_2$ with respect to $x$ at $(t,0)$. Then by \eqref{eq:C2alpha-Laplacian} and \eqref{eq:Domega-estimates}, for $\kappa \in (0,1/2)$, we have 
	\begin{equation}\label{eq:v2-nondiv-estimate}
		\norm{v_2-\hat{q}_2}{\Leb{\infty}\Leb{2}(Q_{\kappa r})}\leq N\kappa^{2+\alpha+d/2}\norm{v_2}{\Leb{\infty}\Leb{2}(Q_r)}.
	\end{equation}
	Then following exactly the same argument as in (i) using \eqref{eq:v1-estimate-nondiv}, \eqref{eq:v2-nondiv-estimate}, and Lemma \ref{lem:polynomial-comparison} (ii) instead, we have 
	\begin{align*}
		\phi_2(\kappa r)\leq N_0 \kappa^{\alpha} \phi_2(r)+ N (\rho_A(r)\norm{D^2u}{\Leb{\infty}(Q_r)} + \rho_f(r)+\widehat{\rho}_{Dg}(r)+\psi(r))
	\end{align*} 
	for some constants $N_0=N_0(d,\nu)>0$ and $N=N(d,\nu,\kappa)>0$.
	
	Then by repeating the same line of the proof of Proposition \ref{prop:velocity-approximation} (i) using \eqref{eq:D-vorticity-Dini-estimate} instead, for any $X_0=(t_0,x_0)\in Q_{3/4}$, we have
	\begin{align*}
		\phi_2(\kappa^j r,X_0)&\apprle_{d,\nu,\beta} \kappa^{j\beta}\left(\phi_2(r,X_0)+\psi(r,X_0)\right)+\norm{D^2u}{\Leb{\infty}(Q_r(X_0))}I_{\Rho_A}(r)\\
		&\relphantom{=}+I_{\Rho_f}(r)+I_{\widehat{\Rho}_{Dg}}(r).
	\end{align*}
	This completes the proof of Proposition \ref{prop:velocity-approximation} (ii).\hfill\qed
	
	\begin{remark}
		Following a similar argument as in Remark \ref{rem:rate-of-psi}, we get the following estimates: for $\beta \in (0,1)$ and $0<s\leq r$, we have
		\begin{equation*}
			\begin{aligned}
				\phi_2(s,X_0)&\apprle \left(\frac{s}{r}\right)^\beta \left( \phi_2(r,X_0)+ \psi(r,X_0) \right)\\
				&\relphantom{=}+\norm{D^2 u}{\Leb{\infty}(Q_r(X_0))}I_{\Rho_A}(s)+I_{\Rho_f}(s)+I_{\widehat{\Rho}_{Dg}}(s).
			\end{aligned}
		\end{equation*}
	\end{remark}
	
	\section{Proof of Theorems \ref{thm:A} and \ref{thm:B}}\label{sec:proof-main}
	
	This section is devoted to proving Theorems \ref{thm:A} and \ref{thm:B}.
	
	\begin{proof}[Proof of Theorem \ref{thm:A}]
		By approximation argument, we may assume that $u \in C^0_tC^1(\overline{Q_1})$, and $a^{ij}$, $\boldF$, and $g$ are smooth. This reduction will be justified in Appendix \ref{app:approximation}. Since $u \in C^0_t C^1(\overline{Q_1})$, it follows that 
		\[ \phi_1(r,t_0,x_0)=\sup_{t\in (t_0-r^2,t_0]} \inf_{p\in\mathcal{P}_1} \frac{c}{r}\left(\fint_{B_r(x_0)} |u(t,x)-p(x)|^2\myd{x}\right)^{1/2}\]
		for some constant $c=c(d)>0$.
		
		For $X_0=(t_0,x_0) \in Q_{3/4}$ and $r\in (0,1/4)$, choose a polynomial $q_{(t,x_0),r} \in \mathcal{P}_1$ so that 
		\[ \phi_1(r,X_0)=\sup_{t\in (t_0-r^2,t_0]} \frac{c}{r}\left(\fint_{B_r(x_0)} |u(t,x)-q_{(t,x_0),r}(x)|^2 \myd{x}\right)^{1/2}.\]
		
		Write $q_{(t,x_0),r}(x)=a_{(t,x_0),r}+b_{(t,x_0),r}\cdot (x-x_0)$. Then 
		\begin{equation}\label{eq:u-q-average}
			|q_{(t,x_0),r}(x)-q_{(t,x_0),\kappa r}(x)|^2\leq 2|u(t,x)-q_{(t,x_0),r}(x)|^2 + 2|u(t,x)-q_{(t,x_0),\kappa r}(x)|^2.
		\end{equation}
		By taking average of \eqref{eq:u-q-average} over $x\in B_{\kappa r}(x_0)$, we get 
		\begin{align*}
			&\left(\fint_{B_{\kappa r}(x_0)} |q_{(t,x_0),r}(x)-q_{(t,x_0),\kappa r}(x)|^2 \myd{x}\right)^{1/2}\\
			&\apprle \kappa^{-d/2} \left(\fint_{B_r(x_0)} |u(t,x)-q_{(t,x_0),r}(x)|^2 \myd{x}\right)^{1/2}\\
			&\relphantom{=} + \left(\fint_{B_{\kappa r}(x_0)} |u(t,x)-q_{(t,x_0),\kappa r}(x)|^2 \myd{x}\right)^{1/2}
		\end{align*}
		for $t\in (t_0-(\kappa r)^2,t_0]$. This implies that 
		\begin{equation}\label{eq:polynomial-average}
			\begin{aligned}
				&\frac{1}{\kappa r}\left(\fint_{B_{\kappa r}(x_0)} |q_{(t,x_0),r}(x)-q_{(t,x_0),\kappa r}(x)|^2 \myd{x}\right)^{1/2}\\
				&\apprle \kappa^{-d/2-1} \phi_1(r,X_0)+\phi_1(\kappa r,X_0)
			\end{aligned}
		\end{equation}
		for $t\in (t_0-(\kappa r)^2,t_0]$. 
		
		Since 
		\[ q_{(t,x_0),r}(x)-q_{(t,x_0),\kappa r}(x)=(a_{(t,x_0),r}-a_{(t,x_0),\kappa r})+(b_{(t,x_0),r}-b_{(t,x_0),\kappa r})\cdot (x-x_0),\]
		a change of variable gives 
		\begin{align*}
			&\frac{1}{\kappa r} \left(\fint_{B_{\kappa r}(x_0)} |q_{(t,x_0),r}(x)-q_{(t,x_0),\kappa r}(x)|^2\myd{x} \right)^{1/2}\\
			&=\frac{1}{\kappa r} \left(\fint_{B_{\kappa r}(x_0)} |a_{(t,x_0),r}-a_{(t,x_0),\kappa r}|^2+|b_{(t,x_0),r}-b_{(t,x_0),\kappa r}|^2 |x-x_0|^2\myd{x} \right)^{1/2}\\
			&\apprge |b_{(t,x_0),r}-b_{(t,x_0),\kappa r}|
		\end{align*}
		for $t\in (t_0-(\kappa r)^2,t_0]$. Hence by \eqref{eq:polynomial-average}, for $t\in (t_0-(\kappa r)^2,t_0]$, we have
		\begin{equation*}\label{eq:vector-field-b-iteration}
			|b_{(t,x_0),r}-b_{(t,x_0),\kappa r}|\apprle \kappa^{-d/2-1} \phi_1(r,X_0)+\phi_1(\kappa r,X_0),
		\end{equation*}
		where the implicit constant depends only on $d$. 
		
		By iteration and the triangle inequality, for $j\in \mathbb{N}$, we have
		\begin{equation}\label{eq:b-x0-r}
			|b_{(t,x_0),\kappa^j r}-b_{(t,x_0),r}|\apprle (\kappa^{-d/2-1} +1)\sum_{l=0}^j \phi_1(\kappa^l r,X_0)
		\end{equation}
		for $t\in (t_0-(\kappa^j r)^2,t_0]$.
		
		By a priori assumption $u\in C^0_tC^1(\overline{Q_1})$, we have 
		\begin{equation}\label{eq:limit-Du-X0}
			\lim_{\rho \rightarrow 0+} \phi_1(\rho,X_0)=0\quad \text{and}\quad \lim_{j\rightarrow\infty} b_{(t_0,x_0),\kappa^j r}=Du(X_0).
		\end{equation}
		Moreover, it follows from  \eqref{eq:phi-psi-iteration-final}, \eqref{eq:b-x0-r}, and \eqref{eq:limit-Du-X0} that 
		\begin{equation}\label{eq:Du-b-difference}
			\begin{aligned}
				&|Du(X_0)-b_{(t_0,x_0),r}|\\
				&\apprle (\kappa^{-d/2}+1)\sum_{l=0}^\infty \phi_1(\kappa^l r,X_0)\\
				&\apprle \phi_1(r,X_0)+\psi(\omega,r,X_0)+\norm{Du}{\Leb{\infty}(Q_r(X_0))}I_{\Rho_A}(r)+I_{\Rho_{\boldF}}(r)+I_{\widehat{\Rho}_g}(r).
			\end{aligned}
		\end{equation}
		On the other hand, for $t\in (t_0-r^2,t_0]$, we note that 
		\begin{align}
			|b_{(t,x_0),r}|
			&\apprle \frac{1}{r}\left(\fint_{B_r(x_0)} |q_{(t,x_0),r}|^2 \myd{x}\right)^{1/2}\label{eq:b-X0-estimate}\\
			&\apprle \frac{1}{r}\left(\fint_{B_r(x_0)} |u(t,x)-q_{(t,x_0),r}(x)|^2 \myd{x}\right)^{1/2}+\frac{1}{r}\left(\fint_{B_r(x_0)} |u(t,x)|^2 \myd{x}\right)^{1/2}\nonumber\\
			&\apprle \phi_1(r,X_0)+r^{-1-d/2} \norm{u}{\Leb{\infty}\Leb{2}(Q_r(X_0))}.\nonumber
		\end{align}
		
		Since 
		\[ \phi_1(r,X_0)+\psi(\omega,r,X_0)\apprle r^{-1-d/2}\left(\norm{u}{\Leb{\infty}\Leb{2}(Q_r(X_0))}+\norm{Du}{\Leb{2}(Q_r(X_0))}\right),
		\]
		it follows from \eqref{eq:Du-b-difference} and \eqref{eq:b-X0-estimate} that
		\begin{align*}
			|Du(X_0)|&\apprle r^{-1-d/2}\left(\norm{u}{\Leb{\infty}\Leb{2}(Q_r(X_0))}+\norm{Du}{\Leb{2}(Q_r(X_0))}\right)\\
			&\relphantom{=} +\norm{Du}{\Leb{\infty}(Q_r(X_0))}I_{\Rho_A}(r)+ I_{\Rho_{\boldF}}(r)+I_{\widehat{\Rho}_g}(r).
		\end{align*}

			By taking the supremum for $X_0 \in Q_r(X_1)$, where $X_1\in Q_{1/4}$ and $r\in (0,1/4)$, we have 
			\begin{equation}\label{eq:Du-estimate-preparation}
				\begin{aligned}
					\norm{Du}{\Leb{\infty}(Q_r(X_1))}&\leq N_2 r^{-1-d/2}\left(\norm{u}{\Leb{\infty}\Leb{2}(Q_{2r}(X_1))}+\norm{Du}{\Leb{2}(Q_{2r}(X_1))}\right)\\
					&\relphantom{=}+N_2\norm{Du}{\Leb{\infty}(Q_{2r}(X_1))}I_{\Rho_A}(r)+N_2(I_{\Rho_{\boldF}}(r)+I_{\widehat{\Rho}_g}(r))
				\end{aligned}
			\end{equation}
			for some constant $N_2=N_2(d,\nu)>0$.
			
			Choose $r_0 \in (0,1/4)$ so that 
			\[  N_2 I_{\Rho_A}(r_0)\leq 3^{-d/2-1}. \]
			Then for $X_1 \in Q_{1/4}$ and $r\in (0,r_0]$, \eqref{eq:Du-estimate-preparation} gives
			\begin{align*}
				\norm{Du}{\Leb{\infty}(Q_r(X_1))}&\leq 3^{-(d/2+1)}\norm{Du}{\Leb{\infty}(Q_{2r}(X_1))}+N(I_{\Rho_{\boldF}}(r_0)+I_{\widehat{\Rho}_g}(r_0))\\
				&\relphantom{=}+Nr^{-1-d/2}\left(\norm{u}{\Leb{\infty}\Leb{2}(Q_{2r}(X_1))}+\norm{Du}{\Leb{2}(Q_{2r}(X_1))}\right)
			\end{align*}
			for some constant $N=N(d,\nu)>0$.
			
			For $k=1,2,\ldots$, define $r_k = 3/4-2^{-k}$. Note that $r_{k+1}-r_k=2^{-(k+1)}$ for $k\geq 1$ and $r_1=1/4$. For $X_1 \in Q_{r_k}$ and $r= 2^{-k-2}$, we have $Q_{2r}(X_1)\subset Q_{r_{k+1}}$. We take $k_0\geq 1$ sufficiently large so that $2^{-k_0-2}\leq r_0$. It then follows that for any $k\geq k_0$, 
			\begin{align*}
				\norm{Du}{\Leb{\infty}(Q_{r_k})}&\leq 3^{-(d/2+1)}\norm{Du}{\Leb{\infty}(Q_{r_{k+1}})}+N(I_{\Rho_{\boldF}}(r_0)+I_{\widehat{\Rho}_g}(r_0))\\
				&\relphantom{=}+N2^{(d/2+1)k}\left(\norm{u}{\Leb{\infty}\Leb{2}(Q_{r_{k+1}})}+\norm{Du}{\Leb{2}(Q_{r_{k+1}})}\right).
			\end{align*}
			By multiplying the above by $3^{-(d/2+1)k}$ and then summing the terms with respect to $k=k_0$, $k_0+1$, ..., we have
			\begin{align*}
				&\sum_{k=k_0}^\infty 3^{-(d/2+1)k} \norm{Du}{\Leb{\infty}(Q_{r_k})}\leq \sum_{k=k_0}^\infty 3^{-(d/2+1)(k+1)}\norm{Du}{\Leb{\infty}(Q_{r_{k+1}})}\\
				&\relphantom{=}+N\norm{u}{\Leb{\infty}\Leb{2}(Q_{3/4})}+N\norm{Du}{\Leb{2}(Q_{3/4})}+N(I_{\Rho_{\boldF}}(1)+I_{\widehat{\Rho}_g}(1))
			\end{align*}
			for some constant $N=N(d,\nu)>0$.  Since $u\in C^0_tC^1(Q_1)$, the summations of the above inequality are finite. Hence it follows that 
			\begin{equation}\label{eq:boundedness-gradient-estimate}
				\norm{Du}{\Leb{\infty}(Q_{1/4})}\apprle \norm{u}{\Leb{\infty}\Leb{2}(Q_{3/4})}+\norm{Du}{\Leb{2}(Q_{3/4})}+I_{\Rho_{\boldF}}(1)+I_{\widehat{\Rho}_g}(1).
			\end{equation}
			Moreover, it follows from the Caccioppoli type estimates for Stokes equations with variable coefficients (see Dong-Phan \cite[Theorem 1.8]{DP21}) that 
			\begin{equation}\label{eq:Du-estimate-Caccioppoli}
				\begin{aligned}
					\norm{Du}{\Leb{2}(Q_{3/4})}\apprle \norm{u}{\Leb{\infty}\Leb{2}(Q_1)}+\norm{\boldF}{\Leb{2}(Q_1)}+\norm{g}{\Leb{2}(Q_1)}.
				\end{aligned}
			\end{equation}
			Hence by \eqref{eq:boundedness-gradient-estimate} and \eqref{eq:Du-estimate-Caccioppoli}, we have 
			\begin{equation}\label{eq:boundedness-gradient-estimate-2}
				\begin{aligned}
					\norm{Du}{\Leb{\infty}(Q_{1/4})}&\apprle \norm{u}{\Leb{\infty}\Leb{2}(Q_{1})}+\norm{\boldF}{\Leb{2}(Q_{1})}+\norm{g}{\Leb{2}(Q_1)}+I_{\Rho_{\boldF}}(1)+I_{\widehat{\Rho}_g}(1).
				\end{aligned}
			\end{equation}

			Next, we show the modulus of continuity of $Du$ with respect to $x$. By \eqref{eq:Du-estimate-preparation}, for $t\in (-1/64,0)$ and $x,y \in B_{1/8}$ with $r:=|x-y|\in (0,1/8)$, we have
			\begin{equation}\label{eq:Du-average}
				\begin{aligned}
					&|Du(t,x)-Du(t,y)|\\
					&\leq|Du(t,x)-b_{(t,x),r}| + |b_{(t,x),r}-b_{(t,y),r}| + |b_{(t,y),r}-Du(t,y)|.
				\end{aligned}
			\end{equation}
			
			If we write $\overline{x}=(x+y)/2$, then since $B_{r/2}(\overline{x})\subset B_{r}(x)\cap B_r(y)$, it follows that 
			\begin{align*}
				|b_{(t,x),r}-b_{(t,y),r}|&\apprle \frac{1}{r}\left(\fint_{B_{r/2}(\overline{x})} |q_{(t,x),r}-q_{(t,y),r}|^2 \myd{x}dt\right)^{1/2} \\
				&\apprle \frac{1}{r} \left(\fint_{B_r(x)\cap B_r(y)} |q_{(t,x),r}-q_{(t,y),r}|^2 \myd{x}dt\right)^{1/2}\\
				&\apprle \phi_1(r,t,x)+\phi_1(r,t,y).
			\end{align*}
			Hence by \eqref{eq:Du-b-difference} and \eqref{eq:Du-average}, we have 
			\begin{align}\label{eq:Du-estimates}
				|Du(t,x)-Du(t,y)|&\apprle \sup_{x_0 \in B_{1/8}} |Du(t,x_0)-b_{(t,x_0),r}|+\phi_1(r,t,x)+\phi_1(r,t,y)\nonumber\\
				&\apprle \sup_{x_0\in B_{1/8}} \sum_{j=0}^\infty \phi_1(2^{-j}r,t,x_0)\nonumber\\
				&\apprle \sup_{x_0 \in B_{1/8}} \left(\phi_1(r,t,x_0)+\psi(r,t,x_0)\right)\\
				&\relphantom{=}+\norm{Du}{\Leb{\infty}(Q_{1/4})}I_{\Rho_A}(r)+I_{\Rho_{\boldF}}(r)+I_{\widehat{\Rho}_g}(r)\nonumber\\
				&\apprle \sup_{x_0 \in B_{1/8}} r^\beta (\phi_1(1/2,t,x_0)+\psi(1/2,t,x_0))\nonumber\\
				&\relphantom{=}+\norm{Du}{\Leb{\infty}(Q_{1/4})}I_{\Rho_A}(r)+I_{\Rho_{\boldF}}(r)+I_{\widehat{\Rho}_g}(r),\nonumber
			\end{align}
			where we used \eqref{eq:phi-psi-iteration-final} in the last inequality.
			
			Hence by \eqref{eq:phi-estimate}, \eqref{eq:boundedness-gradient-estimate-2}, \eqref{eq:Du-estimates}, and Dong-Phan \cite[Theorem 1.8]{DP21}, we have 
			\begin{equation}\label{eq:Du-uniformly-conti}
				\begin{aligned}
					&|Du(t,x)-Du(t,y)|\\
					&\apprle |x-y|^\beta \left(\norm{u}{\Leb{\infty}\Leb{2}(Q_1)}+\norm{\boldF}{\Leb{2}(Q_1)}+\norm{g}{\Leb{2}(Q_1)}\right)+I_{\Rho_{\boldF}}(|x-y|)+I_{\widehat{\Rho}_{g}}(|x-y|)\\
					&\relphantom{=}+\left(\norm{u}{\Leb{\infty}\Leb{2}(Q_1)}+\norm{\boldF}{\Leb{2}(Q_1)}+\norm{g}{\Leb{2}(Q_1)}+I_{\Rho_{\boldF}}(1)+I_{\widehat{\Rho}_g}(1) \right)I_{\Rho_A}(|x-y|)
				\end{aligned}
			\end{equation}
			for any $x,y\in B_{1/8}$ with $|x-y|<1/8$ and $t\in (-1/64,0)$. This completes the proof of Theorem \ref{thm:A}. 
		\end{proof}
		
		{\begin{remark}\label{rem:continuity-vorticity}
		Following the argument as in the proof of Theorem \ref{thm:A} (or \cite{DEK21}), it follows from \eqref{eq:rate-vorticity} that the vorticity $\omega$ satisfies
		\begin{align*}
		&|\omega(t,x)-\omega(s,y)|\\
		&\apprle R^\beta (\norm{u}{\Leb{\infty}\Leb{2}(Q_1)}+\norm{\boldF}{\Leb{2}(Q_1)}+\norm{g}{\Leb{2}(Q_1)})\\
		&\relphantom{=}+(\norm{u}{\Leb{\infty}\Leb{2}(Q_1)}+\norm{\boldF}{\Leb{2}(Q_1)}+\norm{g}{\Leb{2}(Q_1)}+I_{\widehat{\Rho}_{\boldF}}(1)+I_{\widehat{\Rho}_g}(1) )I_{\Rho_A}(R)+I_{\Rho_{g}}(R),
		\end{align*}
		where $R=|t-s|^{1/2}+|x-y|$ and $(t,x),(s,y)\in Q_{1/8}$ satisfying $R<1/8$.
		
		 Similarly, one can show that $D\omega$ is uniformly continuous in $(t,x)$ for the equation in nondivergence form.
		\end{remark}
		}
		
		\begin{remark}\label{rem:divergence-Holder} 
			Suppose that $A\in C^\alpha_x$, $\boldF,g\in C^\alpha_x(Q_2)$ for some $\alpha \in (0,1)$. Then by \eqref{eq:Du-uniformly-conti} with $\beta=\alpha$, we have 
			\begin{equation}\label{eq:Du-Holder-conti}
				\begin{aligned}
					&|Du(t,x)-Du(t,y)|\\
					&\apprle (\norm{u}{\Leb{\infty}\Leb{2}(Q_1)}+\norm{\boldF}{C^\alpha_x(Q_1)}+\norm{g}{C^\alpha_x(Q_1)})(1+[A]_{C^\alpha_x})|x-y|^\alpha
				\end{aligned}
			\end{equation}
			for any $t\in (-1/64,0)$ and $x,y\in B_{1/8}$. This shows that 
			\begin{align*}
				[Du]_{C^\alpha_x(Q_{1/8})}&\apprle \norm{u}{\Leb{\infty}\Leb{2}(Q_1)}+\norm{\boldF}{C^\alpha_x(Q_1)}+\norm{g}{C^\alpha_x(Q_1)},
			\end{align*}
			where the implicit constant depends on $d$, $\nu$, $\alpha$, and $[A]_{C^\alpha_x}$. 
			Similarly, it follows from \eqref{eq:rate-vorticity} that $\omega \in C^{\alpha/2,\alpha}_{t,x}(Q_{1/8})$ and 
			\begin{align*}
				[\omega]_{C^{\alpha/2,\alpha}_{t,x}(Q_{1/8})}&\apprle \norm{u}{\Leb{\infty}\Leb{2}(Q_1)}+\norm{\boldF}{C^\alpha_x(Q_1)}+\norm{g}{C^\alpha_x(Q_1)}.
			\end{align*}
		\end{remark}
		\begin{remark}\label{rem:further-regularity}
			Suppose in addition that $u\in C^{\alpha/2,\alpha}_{t,x}(Q_1)$ is a solution of \eqref{eq:div-form} with $A\in C^{\alpha/2,\alpha}_{t,x}$ and $\boldF,g \in C^{\alpha/2,\alpha}_{t,x}(Q_2)$ for some $\alpha \in (0,1)$. Then $Du\in C^{\alpha^-/2,\alpha}_{t,x}(Q_{1/8})$. 
			
			For $\beta \in (0,\alpha/2)$ and $h\in (0,1/8)$, define 
			\[ \delta_h u(t,x)=\frac{u(t,x)-u(t-h,x)}{h^{\beta}}\quad \text{and}\quad \tau_h u(t,x)=u(t-h,x).\]
			Then it is easy to check that for $0<h<r$, we have 
			\[ [\delta_h u]_{C^{\alpha/2-\beta,\alpha-2\beta}_{t,x}(Q_{r-h})}\leq 2[u]_{C^{\alpha/2,\alpha}_{t,x}(Q_r)}\]
			and $\delta_h u$ satisfies
			\begin{equation*}
				\left\{
				\begin{aligned}
					\partial_t (\delta_h u)-D_i(a^{ij}D_j(\delta_h u))+\nabla (\delta_h \pi)&=D_i((\delta_h a^{ij})D_j(\tau_h u))+\Div(\delta_h \boldF),\\
					\Div (\delta_h u)&=\delta_h g.
				\end{aligned}
				\right.
			\end{equation*}
			Hence for sufficiently small $h$, it follows from \eqref{eq:boundedness-gradient-estimate-2} that 
			\begin{equation}\label{eq:Du-delta-h}
				\begin{aligned}
					\norm{D(\delta_h u)}{\Leb{\infty}(Q_{1/4})}&\apprle \norm{\delta_h u}{\Leb{\infty}\Leb{2}(Q_{7/8})}+\norm{\delta_h \boldF}{\Leb{2}(Q_{7/8})}+\norm{\delta_h g}{\Leb{2}(Q_{7/8})}\\
					&\relphantom{=}+[\delta_h\boldF]_{C^{\alpha-2\beta}_x(Q_{7/8})}+[\delta_h g]_{C^{\alpha-2\beta}_x(Q_{7/8})}\\
					&\apprle [u]_{C^{\alpha/2,\alpha}_{t,x}(Q_1)}+[\boldF]_{C^{\alpha/2,\alpha}_{t,x}(Q_2)}+[g]_{C^{\alpha/2,\alpha}_{t,x}(Q_2)}\\
					&\relphantom{=} +[\delta_h\boldF]_{C^{\alpha-2\beta}_x(Q_{7/8})}+[\delta_h g]_{C^{\alpha-2\beta}_x(Q_{7/8})}\\
					&\apprle [u]_{C^{\alpha/2,\alpha}_{t,x}(Q_1)}+[\boldF]_{C^{\alpha/2,\alpha}_{t,x}(Q_2)}+[g]_{C^{\alpha/2,\alpha}_{t,x}(Q_2)},
				\end{aligned}
			\end{equation}
			where the implicit constant is independent of $h$. 
			
			By \eqref{eq:Du-Holder-conti} and \eqref{eq:Du-delta-h}, we have 
			\begin{align*}
				&|Du(t,x)-Du(t-h,y)|\\
				&\leq |Du(t,x)-Du(t,y)|+|Du(t-h,y)-Du(t,y)|\\
				&\apprle (\norm{u}{\Leb{\infty}\Leb{2}(Q_1)}+\norm{\boldF}{\Leb{2}(Q_1)}+\norm{g}{\Leb{2}(Q_1)})(1+[A]_{C^\alpha})|x-y|^\alpha\\
				&\relphantom{=}+\left([u]_{C^{\alpha/2,\alpha}_{t,x}(Q_1)}+\norm{\boldF}{C^{\alpha/2,\alpha}_{t,x}(Q_2)}+\norm{g}{C^{\alpha/2,\alpha}_{t,x}(Q_2)} \right)h^\beta,
			\end{align*}
			which proves that $Du \in C^{\alpha^{-}/2,\alpha}_{t,x}(Q_{1/8})$. 
		\end{remark}

		\begin{proof}[Proof of Theorem \ref{thm:B}]
			Following a similar argument as in the proof of Theorem \ref{thm:A}, we may assume that $u\in C^0_tC^2(\overline{Q_1})$ and $a^{ij}$, $f$, and $g$ are smooth. Also, following the same line as in Theorem \ref{thm:A} using Proposition \ref{prop:velocity-approximation} (ii) instead of (i), we get 
			\begin{equation}\label{eq:D2u-estimate-apriori}
				\norm{D^2u}{\Leb{\infty}(Q_{1/4})}\apprle \norm{u}{\Leb{\infty}\Leb{2}(Q_{3/4})}+ \norm{u}{\Sob{0,2}{2}(Q_{3/4})}+I_{\Rho_f}(1)+I_{\widehat{\Rho}_{Dg}}(1).
			\end{equation}
			Moreover, it follows from Dong-Phan \cite[Theorem 1.10]{DP21} that 
			\begin{equation}\label{eq:D2u-Caccioppoli-inequality-estimate}
				\norm{u}{\Sob{0,2}{2}(Q_{3/4})}\apprle \norm{u}{\Leb{\infty}\Leb{2}(Q_1)}+\norm{f}{\Leb{2}(Q_1)}+\norm{Dg}{\Leb{2}(Q_1)}.
			\end{equation}
			Hence, by \eqref{eq:D2u-estimate-apriori} and \eqref{eq:D2u-Caccioppoli-inequality-estimate}, we have 
			\begin{equation}\label{eq:D2u-estimate-apriori-2}
				\begin{aligned}
					\norm{D^2u}{\Leb{\infty}(Q_{1/4})}\apprle \norm{u}{\Leb{\infty}\Leb{2}(Q_{1})}+\norm{f}{\Leb{2}(Q_1)}+\norm{Dg}{\Leb{2}(Q_1)}+I_{\Rho_f}(1)+I_{\widehat{\Rho}_{Dg}}(1).
				\end{aligned}
			\end{equation}
			
			We also have the following estimate for the modulus of the continuity of $D^2 u$:  
			fix $\beta \in (0,1)$. Then for $t\in (-1/64,0)$, $x,y\in B_{1/8}$ with $|x-y|<1/8$, we have
			\begin{align}\label{eq:Du-uniformly-conti-nondiv}
				&|D^2 u(t,x)-D^2u(t,y)|\nonumber\\
				&\apprle |x-y|^\beta\left( \norm{u}{\Leb{\infty}\Leb{2}(Q_1)}+\norm{f}{\Leb{2}(Q_1)}+\norm{Dg}{\Leb{2}(Q_1)}\right)\\
				&\relphantom{=}+I_{\Rho_f}(|x-y|)+I_{\widehat{\Rho}_{Dg}}(|x-y|)\nonumber\\
				&\relphantom{=}+\left(\norm{u}{\Leb{\infty}\Leb{2}(Q_1)}+\norm{f}{\Leb{2}(Q_1)}+\norm{Dg}{\Leb{2}(Q_1)}+I_{\Rho_f}(1)+I_{\widehat{\Rho}_{Dg}}(1) \right)I_{\Rho_A}(|x-y|).\nonumber
			\end{align}
			This completes the proof of Theorem \ref{thm:B}.
		\end{proof}

		\begin{remark}\label{rem:nondivergence-Holder}\leavevmode 
			Suppose that $A\in C^\alpha_x$ and $f,Dg \in C^\alpha_x(Q_2)$ for some $\alpha \in (0,1)$. Then by \eqref{eq:Du-uniformly-conti-nondiv}, we have
			\begin{align*}
				&|D^2u(t,x)-D^2u(t,y)|\\
				&\apprle (\norm{u}{\Leb{\infty}\Leb{2}(Q_1)}+\norm{f}{C^\alpha_x(Q_1)}+\norm{Dg}{C^\alpha_x(Q_1)})(1+[A]_{C^\alpha_x})|x-y|^\alpha
			\end{align*}
			for any $t\in (-1/64,0)$ and $x,y\in B_{1/4}$. This implies that 
			\begin{align*}
				[D^2u]_{C^\alpha_x(Q_{1/8})}&\apprle \norm{u}{\Leb{\infty}\Leb{2}(Q_1)}+\norm{f}{C^\alpha_x(Q_2)}+\norm{Dg}{C^\alpha_x(Q_2)},
			\end{align*}
			where the implicit constant depends on $d$, $\nu$, $\alpha$, and $[A]_{C^\alpha_x}$. Similarly, it follows from \eqref{eq:rate-D-vorticity} that $D\omega \in C^{\alpha/2,\alpha}_{t,x}(Q_{1/8})$ and 
			\begin{align*}
				[D\omega]_{C^{\alpha/2,\alpha}_{t,x}(Q_{1/8})}&\apprle \norm{u}{\Leb{\infty}\Leb{2}(Q_1)}+\norm{f}{C^\alpha_x(Q_1)}+\norm{Dg}{C^\alpha_x(Q_1)}.
			\end{align*}
		\end{remark}
		
		\begin{remark}\label{rem:nondivergence-Holder-additional}\leavevmode 
			Suppose that $A\in C^{\alpha/2,\alpha}_{t,x}$ and $f,Dg \in C^{\alpha/2,\alpha}_{t,x}(Q_2)$ for some $\alpha \in (0,1)$. Following a similar argument as in Remark \ref{rem:further-regularity}, one can show that if $u\in C^{\alpha/2,\alpha}_{t,x}(Q_1)$ in addition, then $D^2 u\in C^{\alpha^-/2,\alpha}_{t,x}(Q_{1/8})$.
		\end{remark}
		
		\section*{Acknowledgment}
		H. Dong would like to thank Prof. Ping Zhang for valuable discussions on potential applications of the results from years ago. 
		
		\subsection*{Conflict of interest}
		
		The authors declare that there are no conflict of interest to this paper.

		\subsection*{Data availibility}
		The authors declare that data sharing is not applicable to this paper since there is no datasets were generated during the current study.
		\appendix
		
		\section{Solvability of Stokes equations with simple coefficients}\label{sec:Appendix}
		In this section, we prove solvability results for Stokes equations with simple coefficients. Such result will be used in the proof of Proposition \ref{prop:velocity-approximation}. 
		
		The following lemma will be used to estimate a solution $u$ of \eqref{eq:div-form} in $\Leb{\infty}\Leb{2}(Q_r)$. This lemma holds for any bounded viscosity coefficients which are uniformly elliptic.
		
		\begin{lemma}\label{lem:whole-space-div-free}
			Let $u\in \mathcal{H}^1_2((-r^2,0)\times\mathbb{R}^d)$ be a weak solution of 
			\begin{equation}\label{eq:weak-sol-div-free}
				\partial_t u-D_i(a^{ij}D_j u)+\nabla \pi =f+\Div \boldF\quad \text{and}\quad \Div u=0\quad \text{in } (-r^2,0)\times \mathbb{R}^d 
			\end{equation}
			for some function $f,\boldF \in \Leb{2}((-r^2,0)\times \mathbb{R}^d)$ with {$\supp f \subset (-r^2,0)\times B_{2r}$} and $u(-r^2,\cdot)=0$. Then $u \in \Leb{\infty}\Leb{2}((-r^2,0)\times\mathbb{R}^d)$. Moreover, we have 
			\[ \norm{u}{\Leb{\infty}\Leb{2}((-r^2,0)\times\mathbb{R}^d)} +\norm{Du}{\Leb{2}((-r^2,0)\times\mathbb{R}^d)}\leq N \left(r\norm{f}{\Leb{2}((-r^2,0)\times B_{2r})}+\norm{\boldF}{\Leb{2}((-r^2,0)\times\mathbb{R}^d)}\right) \]
			for some constant $N=N(d,\nu)>0$. 
		\end{lemma}
		\begin{proof}
			By a standard duality argument as in \cite[Theorem 2.2]{S15}, one can show that $u\in \Leb{\infty}\Leb{2}(Q_r)$.

			To obtain the estimate, since $\supp f \subset Q_{2r}$, it follows from  H\"older's inequality that
			\begin{align*}
				\left|\int_{-r^2}^0\int_{\mathbb{R}^d} f\cdot u \myd{x}dt\right|\leq \norm{u}{\Leb{\infty}\Leb{2}(Q_r)}\norm{f}{\Leb{1}\Leb{2}(Q_r)}\leq r\norm{u}{\Leb{\infty}\Leb{2}(Q_r)}\norm{f}{\Leb{2}(Q_r)}.
			\end{align*}
			Since $u$ is a weak solution of \eqref{eq:weak-sol-div-free}, it follows that {
			\begin{equation}\label{eq:weak-formulation}
				\begin{aligned}
					&\int_{-r^2}^t \action{\partial_t u(s,\cdot) ,\phi(s,\cdot)}ds +\int_{-r^2}^t\int_{\mathbb{R}^d} {a}^{ij}D_j u\cdot D_i\phi \myd{x}ds\\
					&=-\int_{-r^2}^t\int_{\mathbb{R}^d}\boldF:\nabla \phi \myd{x}ds+\int_{-r^2}^t\int_{\mathbb{R}^d}f\cdot\phi \myd{x}ds
				\end{aligned}
			\end{equation}
			for all $\phi \in C_0^\infty((-r^2,0)\times\mathbb{R}^d)$ with $\Div \phi =0$ and $t\in (-r^2,0)$.} By density argument, the identity holds for all $\phi \in \Sob{0,1}{2,\sigma}((-r^2,0)\times\mathbb{R}^d)$. Since $u\in \Sob{0,1}{2,\sigma}((-r^2,0)\times\mathbb{R}^d)$, we put $\phi=u$ in \eqref{eq:weak-formulation}. Then we have
			\begin{align*}
				&\frac{1}{2}\norm{u}{\Leb{\infty}\Leb{2}((-r^2,0)\times\mathbb{R}^d)}^2 + \nu \norm{Du}{\Leb{2}((-r^2,0)\times\mathbb{R}^d)}^2\\
				&\leq r\norm{u}{\Leb{\infty}\Leb{2}((-r^2,0)\times\mathbb{R}^d)}\norm{f}{\Leb{2}((-r^2,0)\times B_{2r})}+\norm{\boldF}{\Leb{2}((-r^2,0)\times\mathbb{R}^d)}\norm{Du}{\Leb{2}((-r^2,0)\times\mathbb{R}^d)}.
			\end{align*}
			Hence by Young's inequality, we get 
			\[ 
			\norm{u}{\Leb{\infty}\Leb{2}((-r^2,0)\times\mathbb{R}^d)}+\sqrt{\nu}\norm{Du}{\Leb{2}((-r^2,0)\times\mathbb{R}^d)}\leq N\left(r\norm{f}{\Leb{2}((-r^2,0)\times B_{2r})}+\norm{\boldF}{\Leb{2}((-r^2,0)\times\mathbb{R}^d)} \right)
			\]
			for some constant $N=N(d,\nu)>0$. This completes the proof of Lemma \ref{lem:whole-space-div-free}.
		\end{proof}

		Now we are ready to present the main result in this section. We assume that the viscosity coefficients are simple.

		\begin{lemma}\label{lem:approximate-system-existence}
			Let $r>0$. For any $f,\boldF\in \Leb{2}(Q_r)$ and $g\in \Leb{\infty}\Leb{2}(Q_r)$, there exists a solution $u\in\Leb{\infty}\Leb{2}((-r^2,0)\times\mathbb{R}^d)\cap\Leb{2}\tSob{1}((-r^2,0)\times\mathbb{R}^d)$ satisfying 
			\begin{equation*}
				\left\{
				\begin{alignedat}{2}
					\partial_t u -D_i(a^{ij}D_j u)+\nabla \pi&=\Div \boldF+f&&\quad \text{in } Q_r,\\
					\Div u&=g&&\quad \text{in } Q_r,\\
					u(-r^2,\cdot)&=0&&\quad \text{on } \mathbb{R}^d
				\end{alignedat}
				\right.
			\end{equation*}
			and the estimate
			\begin{align*}
				& \norm{u}{\Leb{\infty}\Leb{2}((-r^2,0)\times \mathbb{R}^d)}+\norm{Du}{\Leb{2}((-r^2,0)\times \mathbb{R}^d)}\\
				&\leq N \left(\norm{\boldF}{\Leb{2}(Q_r)}+r\norm{f}{\Leb{2}(Q_r)}+r\norm{g}{\Leb{\infty}\Leb{2}(Q_r)} \right) 
			\end{align*}
			holds for some constant $N=N(d,\nu)>0$.
		\end{lemma}

		\begin{proof}
			{By scaling argument, we may assume that $r=1$. Since $\boldF\in \Leb{2}(Q_1)$ and $g\in \Leb{\infty}\Leb{2}(Q_1)$, we can extend $\boldF$ and $g$ to be zero outside $Q_1$. Take a mollification of $g$ in $(t,x)$. Then for sufficiently small $\varepsilon$, $g^{(\varepsilon)}$ is compactly supported in $Q_{4/3}$. Following exactly the same argument as in \cite[Lemma 6.1]{DK23}, there exist $h_\varepsilon,H_\varepsilon\in \Leb{2}((-1,0)\times\mathbb{R}^d)$ such that $h_\varepsilon=g^{(\varepsilon)}$ in $(-1,0)\times\mathbb{R}^d$ and 
			\[ -\int_{-1}^0\int_{\mathbb{R}^d} H_\varepsilon^i \cdot \nabla (D_i\psi)\myd{x}dt=\int_{-1}^0\int_{\mathbb{R}^d} h_\varepsilon \partial_t \psi \myd{x}dt\]
			holds for all $\psi \in C_0^\infty([-1,0)\times\mathbb{R}^d)$. 
			
Define 
\[
   \tilde{f}(t,x)=\begin{cases}
     f(t,x)&\quad \text{if } x\in B_1\\
     \alpha(t)&\quad \text{if } x\in B_2\setminus B_1
\end{cases},
\]
where 
\[ \alpha(t)=-\frac{1}{|B_2\setminus B_1|} \int_{B_1} f(t,x)\myd{x}.\]
By construction, we have 
\[ \int_{B_2} \tilde{f}(t,x)\myd{x}=0\quad \text{for } t\in (-1,0).\]
Also, it follows from H\"older's inequality that 
\begin{equation}
 \norm{\tilde{f}(t)}{\Leb{2}(B_2)}\leq \norm{f(t)}{\Leb{2}(B_1)}+\norm{\alpha(t)}{\Leb{2}(B_2\setminus B_1)}\leq \left(1+\frac{|B_1|^{1/2}}{|B_2\setminus B_1|^{1/2}} \right)\norm{f(t)}{\Leb{2}(B_1)}.
\end{equation}
This implies that 
\begin{equation}\label{eq:tilde-f-estimate}
 \norm{\tilde{f}}{\Leb{2}((-1,0)\times B_2)}\leq N\norm{f}{\Leb{2}(Q_1)} 
\end{equation}
for some constant $N=N(d)>0$.
			By \eqref{eq:tilde-f-estimate} and Bogovski\u{\i}'s theorem (see e.g. \cite[Theorem III.3.1]{G00}), there exists $\tilde{\boldF}$ such that $\Div \tilde{\boldF}=\tilde{f}$ in $(-1,0)\times B_2$ and 
			\[ \norm{\tilde{\boldF}}{\Leb{2}((-1,0)\times B_2)}\leq  N\norm{\tilde{f}}{\Leb{2}((-1,0)\times B_2)}\leq N\norm{f}{\Leb{2}(Q_1)}\]	
			for some constant $N=N(d)>0$.
We further extend it to be zero outside $B_2$ and still denote it by $\tilde{\boldF}$. Then by \cite[Theorem 4.3]{DK23}, there exists a unique $v^\varepsilon \in \mathcal{H}^1_2((-1,0)\times \mathbb{R}^d)$ satisfying 
			\begin{equation*}\label{eq:vm-solution-3}
				\left\{
				\begin{alignedat}{2}
					\partial_t v^{\varepsilon}-D_i (a^{ij} D_j v^\varepsilon)+\nabla \pi^\varepsilon &=\Div (\tilde{\boldF}+\boldF)&&\quad \text{in } (-1,0)\times \mathbb{R}^d, \\
					\Div v^\varepsilon &=h_\varepsilon&&\quad \text{in } (-1,0)\times \mathbb{R}^d,\\
					v^\varepsilon(-1,\cdot)&=0&&\quad \text{on } \mathbb{R}^d.
				\end{alignedat}
				\right.
			\end{equation*}
			Moreover, we have
			\begin{equation}\label{eq:Dv-epsilon-estimate}
			\begin{aligned}
				\norm{Dv^\varepsilon}{\Leb{2}((-1,0)\times\mathbb{R}^d)}&\apprle \norm{\tilde{\boldF}}{\Leb{2}((-1,0)\times B_2)}+\norm{\boldF}{\Leb{2}(Q_1)}+\norm{h_\varepsilon}{\Leb{2}((-1,0)\times\mathbb{R}^d)}\\
				&\apprle \norm{f}{\Leb{2}(Q_1)}+\norm{\boldF}{\Leb{2}(Q_1)}+\norm{h_\varepsilon}{\Leb{2}((-1,0)\times\mathbb{R}^d)},
			\end{aligned}
			\end{equation}
			where the implicit constant is independent of $\varepsilon$. 
			
			However, it is unclear whether ${v}^\varepsilon$ is in $\Leb{\infty}\Leb{2}((-1,0)\times\mathbb{R}^d)$. To show this, define 
			\[  \tilde{v}^\varepsilon=v^\varepsilon+\nabla \phi,\quad \text{where } \phi(t,x)=\int_{\mathbb{R}^d} \Gamma(x-y)h_\varepsilon(t,y)\myd{y}\]
			and $\Gamma$ is the fundamental solution of $-\Delta$. 
			Since $h_\varepsilon(t,\cdot)$ is compactly supported in $B_{6}$ for each $t$,  we have $-\Delta \phi=h_\varepsilon$. Moreover, it follows from weak Young's convolution inequality (see e.g. \cite[Proposition 8.9]{F99}) and the Calder\'on-Zygmund estimate that
			\begin{equation}\label{eq:CZ-estimate} 
				\begin{aligned}
					\norm{D \phi}{\Leb{\infty}\Leb{2}((-1,0)\times \mathbb{R}^d)}&\apprle \norm{h_\varepsilon}{\Leb{\infty}\Leb{2}((-1,0)\times \mathbb{R}^d)},\\
					\norm{D^2 \phi}{\Leb{2}((-1,0)\times \mathbb{R}^d)}&\apprle \norm{h_\varepsilon}{\Leb{2}((-1,0)\times \mathbb{R}^d)}.
				\end{aligned}
			\end{equation}
			Also, $\Div\tilde{v}^\varepsilon=0$ in $(-1,0)\times\mathbb{R}^d$. Since $\tilde{v}^\varepsilon$ satisfies
			\begin{equation*}
				\begin{aligned}
					\partial_t\tilde{v}^\varepsilon-D_i(a^{ij}D_j\tilde{v}^\varepsilon)+\nabla {\tilde{\pi}}&=f+\Div\boldF-D_i(a^{ij}D_j \nabla \phi),
				\end{aligned}
			\end{equation*}
			it follows from \eqref{eq:CZ-estimate} and Lemma \ref{lem:whole-space-div-free} that $\tilde{v}^\varepsilon \in \Leb{\infty}\Leb{2}((-1,0)\times\mathbb{R}^d)\cap \Leb{2}\tSob{1}((-1,0)\times\mathbb{R}^d)$ and
			\begin{equation}\label{eq:tilde-v-epsilon-Linfty-estimate}
				\begin{aligned}
					&\norm{\tilde{v}^\varepsilon}{\Leb{\infty}\Leb{2}((-1,0)\times\mathbb{R}^d)}+\norm{D\tilde{v}^\varepsilon}{\Leb{2}((-1,0)\times\mathbb{R}^d)}\\
					&\apprle \norm{f}{\Leb{2}(Q_1)}+ \norm{\boldF}{\Leb{2}(Q_1)}+\norm{h_\varepsilon}{\Leb{2}((-1,0)\times\mathbb{R}^d)},
				\end{aligned}
			\end{equation}
			where the implicit constant depends only on $d$, $\nu$. 
			Since 
			\begin{equation}\label{eq:h-epsilon-control}
				\norm{h_\varepsilon(t)}{\Leb{2}(\mathbb{R}^d)}\apprle \norm{g(t)}{\Leb{2}(\mathbb{R}^d)}\quad \text{for } t\in (-1,0),
			\end{equation}
			it follows from \eqref{eq:CZ-estimate} and \eqref{eq:tilde-v-epsilon-Linfty-estimate} that  
			\begin{equation}\label{eq:v-epsilon-Linfty-L2-control}
				\begin{aligned}
					\norm{v^\varepsilon}{\Leb{\infty}\Leb{2}((-1,0)\times \mathbb{R}^d)}&\apprle \norm{\tilde{v}^\varepsilon}{\Leb{\infty}\Leb{2}((-1,0)\times \mathbb{R}^d)}+\norm{h_\varepsilon}{\Leb{\infty}\Leb{2}((-1,0)\times \mathbb{R}^d)}\\
					& \apprle \norm{f}{\Leb{2}(Q_1)}+\norm{\boldF}{\Leb{2}(Q_1)}+\norm{g}{\Leb{\infty}\Leb{2}(Q_1)},
				\end{aligned}
			\end{equation}
			where the implicit constant depends only on $d$ and $\nu$.
			
			By \eqref{eq:Dv-epsilon-estimate}, \eqref{eq:h-epsilon-control}, and \eqref{eq:v-epsilon-Linfty-L2-control}, we have
			\begin{align*}
				\norm{{v}^\varepsilon}{\Leb{\infty}\Leb{2}((-1,0)\times \mathbb{R}^d)}+\norm{Dv^\varepsilon}{\Leb{2}((-1,0)\times \mathbb{R}^d)}\apprle  \norm{\boldF}{\Leb{2}(Q_1)} +\norm{g}{\Leb{2}(Q_1)}+\norm{g}{\Leb{\infty}\Leb{2}(Q_1)},
			\end{align*}
			where the implicit constant depends only on $d$ and $\nu$. Hence by weak and weak-* compactness, there exists a weak solution $v\in \Leb{\infty}\Leb{2}((-1,0)\times\mathbb{R}^d) \cap \Leb{2}\tSob{1}((-1,0)\times \mathbb{R}^d)$ satisfying 
			\begin{equation*}\label{eq:vm-solution-4}
				\left\{
				\begin{alignedat}{2}
					\partial_t v-D_i (a^{ij}D_j v)+\nabla \pi &=f+\Div \boldF&&\quad \text{in } Q_1, \\
					\Div v &=g&&\quad \text{in } Q_1,\\
					v(-1,\cdot)&=0&&\quad \text{on } \mathbb{R}^d
				\end{alignedat}
				\right.
			\end{equation*}
			which satisfies the following estimate
			\begin{align*} 
				&\relphantom{=}\norm{v}{\Leb{\infty}\Leb{2}((-1,0)\times \mathbb{R}^d)} + \norm{Dv}{\Leb{2}((-1,0)\times \mathbb{R}^d)} \\
				&\leq N \left(\norm{f}{\Leb{2}(Q_1)}+\norm{\boldF}{\Leb{2}(Q_1)}+\norm{g}{\Leb{2}(Q_1)}+\norm{g}{\Leb{\infty}\Leb{2}(Q_1)}\right),
			\end{align*} 
			where $N=N(d,\nu)>0$. This completes the proof of Lemma \ref{lem:approximate-system-existence}. }
		\end{proof}

		\section{Approximation argument}\label{app:approximation}
		This section is dedicated to justifying the a priori assumptions used in the proofs of Theorems \ref{thm:A} and \ref{thm:B}. Let $u^{(\varepsilon)}$ denote the space-time mollification of $u$.
		
		\subsection{Equations in nondivergence form}
		We mollify the equation \eqref{eq:nondiv-form}
		\begin{equation}\label{eq:Stokes-commutator}
			\left\{
			\begin{aligned}
				\partial_t  u^{(\varepsilon)}-a^{ij}_{(\varepsilon)}D_{ij}u^{(\varepsilon)}+\nabla \pi^{(\varepsilon)}&=f^{(\varepsilon)}+[a^{ij}D_{ij}u]^{(\varepsilon)}-a_{(\varepsilon)}^{ij}D_{ij}u^{(\varepsilon)}\\
				\Div u^{(\varepsilon)}&=g^{(\varepsilon)}
			\end{aligned}
			\right.
		\end{equation}
		in $Q_1$. Then by \cite[Theorem 2.5]{DK23}, there exists a unique $u_1^\varepsilon \in \Sob{1,2}{2}((-1,0)\times\mathbb{R}^d)$ satisfying
		\begin{equation*}
			\left\{
			\begin{aligned}
				\partial_t  u_1^\varepsilon-a^{ij}_{(\varepsilon)}D_{ij}u^{\varepsilon}_1+\nabla \pi^{\varepsilon}_1&=h^{\varepsilon}1_{Q_{3/4}}&&\quad \text{in } (-1,0)\times\mathbb{R}^d,\\
				\Div u_1^\varepsilon&=0&&\quad \text{in } (-1,0)\times\mathbb{R}^d,\\
				u_1^\varepsilon(-1,\cdot)&=0&&\quad \text{on } \mathbb{R}^d,
			\end{aligned}
			\right.
		\end{equation*}
		where 
		\[ h^\varepsilon=[a^{ij}D_{ij}u]^{(\varepsilon)}-a_{(\varepsilon)}^{ij}D_{ij}u^{(\varepsilon)}.
		\]
		Moreover, we have 
		\begin{equation}\label{eq:u1-epsilon-estimate}
			\norm{u_1^\varepsilon}{\Sob{1,2}{2}((-1,0)\times\mathbb{R}^d)}\leq N \norm{h^\varepsilon}{\Leb{2}(Q_{3/4})}\end{equation}
		for some constant $N=N(d,\nu,R_0)>0$. Define $u_2^\varepsilon=u^{(\varepsilon)}-u_1^\varepsilon$. Then $u_2^\varepsilon \in \Sob{1,2}{2}(Q_{3/4})$ satisfies
		\begin{equation*}\label{eq:u2-nondiv}
			\left\{
			\begin{aligned}
				\partial_t  u_2^\varepsilon-a^{ij}_{(\varepsilon)}D_{ij}u^{\varepsilon}_2+\nabla \pi^{\varepsilon}_2&=f^{(\varepsilon)}&&\quad \text{in } Q_{3/4},\\
				\Div u_2^\varepsilon&=g^{(\varepsilon)}&&\quad \text{in } Q_{3/4}.
			\end{aligned}
			\right.
		\end{equation*}

		We first show that $D^2 u_2^\varepsilon \in \Leb{2}C^{0,\alpha}(Q_{3/5})$ for some $\alpha \in (0,1)$. For simplicity, we drop $\varepsilon$ in the notation and simply write $u$ instead of $u_2^\varepsilon$.
		
		Since $a$, $f$, and $g$ are smooth, it follows from Dong-Phan \cite[Theorem 1.10]{DP21} that 
		\[ \norm{D^2 u}{\Leb{2}(Q_{r_1})}\apprle \norm{u}{\Leb{2}(Q_{7/8})}+\norm{f}{\Leb{2}(Q_{7/8})}+\norm{Dg}{\Leb{2}(Q_{7/8})}.\]
		Using a method of finite difference, we further get 
		\[ \norm{D^{k+2} u}{\Leb{2}(Q_{r_{k+1}})}\apprle \norm{u}{\Leb{2}(Q_{7/8})}+\norm{f}{\Sob{0,k}{2}(Q_{7/8})}+\norm{g}{\Sob{0,k+1}{2}(Q_{7/8})}\]
		for any $k$, where $3/5<r_{k+1}<r_k<\dots<r_1<7/8$. 
		Hence by the Sobolev embedding theorem, $u \in \Leb{2}C^{2,\alpha}(Q_{3/5})$  for some $\alpha \in (0,1)$. Moreover, we have 
		\begin{equation}\label{eq:Holder-estimate-spatial}
			\norm{u}{\Leb{2}C^{2,\alpha}(Q_{3/5})}\apprle  \norm{u}{\Leb{2}(Q_{7/8})}+\norm{f}{\Sob{0,k}{2}(Q_{7/8})}+\norm{g}{\Sob{0,k+1}{2}(Q_{7/8})}. 
		\end{equation}
		
		For $h> 0$, define $\delta_h u(t,x)=[u(t-h,x)-u(t,x)]/h$. Then 
		\begin{equation*}\label{eq:u2-nondiv-finite-difference}
			\left\{
			\begin{aligned}
				\partial_t  (\delta_h u)-a^{ij}D_{ij}(\delta_h u)+\nabla (\delta_h\pi)&=\delta_h f+(\delta_h a^{ij})D_{ij}u(\cdot +h,\cdot)&&\quad \text{in } Q_{3/5},\\
				\Div (\delta_h u)&=\delta_h g&&\quad \text{in } Q_{3/5}.
			\end{aligned}
			\right.
		\end{equation*}
		Then by \eqref{eq:Holder-estimate-spatial}, we have 
		\[ \norm{\delta_h u}{\Leb{2}C^{2,\alpha}(Q_{4/7})}\apprle \norm{\partial_t u}{\Leb{2}(Q_{7/8})}+\norm{\partial_t f}{\Sob{0,k}{2}(Q_{7/8})}+\norm{\partial_tg}{\Sob{0,k+1}{2}(Q_{7/8})}, \]
		where the implicit constant is independent of $h$. This implies that $\partial_t u\in \Leb{2}C^{2,\alpha}(Q_{1/2})$. Hence $u_2^\varepsilon \in C_t^0C^{2,\alpha}(Q_{1/2})$.
		
		Since 
		\[ \rho_{a_{(\varepsilon)}}(s)\leq \rho_{a}(s),\quad \rho_{f^{(\varepsilon)}}(s)\leq \rho_f (s),\quad \hat{\rho}_{Dg^{(\varepsilon)}}(s)\leq \hat{\rho}_{Dg}(s),\]
		it follows from \eqref{eq:D2u-estimate-apriori-2} that 
		\begin{align*}
			\norm{D^2u_2^\varepsilon}{\Leb{\infty}(Q_{1/4})}&\apprle \norm{u_2^\varepsilon}{\Leb{\infty}\Leb{2}(Q_{3/4})}+\norm{u_2^\varepsilon}{\Sob{0,2}{2}(Q_{3/4})}+I_{\Rho_f}(1)+I_{\widehat{\Rho}_{Dg}}(1)\\
			&\apprle \norm{u}{\Leb{\infty}\Leb{2}(Q_1)}+\norm{f}{\Leb{2}(Q_1)}+\norm{Dg}{\Leb{2}(Q_1)}+\norm{u_1^\varepsilon}{\Sob{1,2}{2}(Q_{3/4})}\\
			&\relphantom{=}+I_{\Rho_f}(1)+I_{\widehat{\Rho}_{Dg}}(1).
		\end{align*}
		By \eqref{eq:u1-epsilon-estimate}, $\norm{u_1^\varepsilon}{\Sob{1,2}{2}(Q_{3/4})}\rightarrow 0$ as $\varepsilon\rightarrow 0+$. Hence $\norm{D^2 u_2^\varepsilon}{\Leb{\infty}(Q_{1/4})}$ is bounded by a constant which is independent of $\varepsilon$. Since $D^2 u_1^\varepsilon \rightarrow 0$ a.e., it follows that $D^2 u_2^\varepsilon\rightarrow D^2 u$ a.e. and hence 
		\[ \norm{D^2 u}{\Leb{\infty}(Q_{1/4})}\apprle \norm{u}{\Leb{\infty}\Leb{2}(Q_{1})}+\norm{u}{\Sob{0,2}{2}(Q_{1})}+I_{\Rho_f}(1)+I_{\widehat{\Rho}_{Dg}}(1).\]
		
		\subsection{Equations in divergence form}
		The proof is similar to the nondivergence form case. Choose a cut-off function $\eta \in C_0^\infty((-(7/8)^2,(7/8)^2)\times B_{7/8})$ so that $\eta=1$ in $Q_{3/4}$. Similar to \eqref{eq:Stokes-commutator}, we mollify the equation \eqref{eq:div-form} to get
		\begin{equation*}\label{eq:Stokes-commutator-div}
			\left\{
			\begin{aligned}
				\partial_t  u^{(\varepsilon)}-D_i(a^{ij}_{(\varepsilon)}D_{j}u^{(\varepsilon)})+\nabla \pi^{(\varepsilon)}&=\Div ({\boldF}^{(\varepsilon)}+\boldH^\varepsilon),\\
				\Div u^{(\varepsilon)}&=g^{(\varepsilon)},
			\end{aligned}
			\right.
		\end{equation*}
		where 
		\[ \boldH^\varepsilon=(H_1^\varepsilon,\dots,H_d^\varepsilon)\quad \text{and}\quad H_i^\varepsilon=(a^{ij}D_ju)^{(\varepsilon)}-a^{ij}_{(\varepsilon)}D_j u^{(\varepsilon)}.\]
		
		Then by \cite[Theorem 2.6]{DK23}, there exists a unique $u_1^\varepsilon \in \mathcal{H}^1_2((-1,0)\times\mathbb{R}^d)$ satisfying
		\begin{equation*}
			\left\{
			\begin{aligned}
				\partial_t  u_1^\varepsilon-D_i(a^{ij}_{(\varepsilon)}D_{j}u^{\varepsilon}_1)+\nabla \pi^{\varepsilon}_1&=\Div(\boldH^\varepsilon \eta)&&\quad \text{in } (-1,0)\times\mathbb{R}^d,\\
				\Div u_1^\varepsilon&=0&&\quad \text{in } (-1,0)\times\mathbb{R}^d,\\
				u_1^\varepsilon(-1,\cdot)&=0&&\quad \text{on } \mathbb{R}^d.
			\end{aligned}
			\right.
		\end{equation*}
		Moreover, we have 
		\begin{equation}\label{eq:u1-epsilon-div-form}
			\norm{u_1^\varepsilon}{\mathcal{H}^1_2((-1,0)\times\mathbb{R}^d)}\leq N\norm{\boldH^\varepsilon}{\Leb{2}(Q_{7/8})}
		\end{equation}
		for some constant $N=N(d,\nu,R_0)>0$. By Lemma \ref{lem:whole-space-div-free}, we have
		\begin{equation}\label{eq:u1-epsilon-Linfty-L2}
			\norm{u_1^\varepsilon}{\Leb{\infty}\Leb{2}((-1,0)\times\mathbb{R}^d)}\leq N \norm{\boldH^\varepsilon}{\Leb{2}(Q_{7/8})}
		\end{equation}
		for some constant $N=N(d,\nu)>0$.

		By \cite[Theorem 2.5]{DK23} with the method of continuity (see \cite[Theorem 4.5]{DK24} for the proof), there exists a unique $\tilde{u}_1^\varepsilon \in \Sob{1,2}{2}((-1,0)\times\mathbb{R}^d)$ satisfying
		\begin{equation*}
			\left\{
			\begin{aligned}
				\partial_t  \tilde{u}_1^\varepsilon-a^{ij}_{(\varepsilon)}D_{ij}\tilde{u}^{\varepsilon}_1-(D_i a^{ij}_{(\varepsilon)})D_j \tilde{u}_1^\varepsilon+\nabla \pi^{\varepsilon}_1&=\Div(\boldH^\varepsilon \eta)&&\quad \text{in } (-1,0)\times\mathbb{R}^d,\\
				\Div \tilde{u}_1^\varepsilon&=0&&\quad \text{in } (-1,0)\times\mathbb{R}^d,\\
				\tilde{u}_1^\varepsilon(-1,\cdot)&=0&&\quad \text{on } \mathbb{R}^d.
			\end{aligned}
			\right.
		\end{equation*}
		Hence by the uniqueness result in \cite[Theorem 2.6]{DK23}, $u_1^\varepsilon \in \Sob{1,2}{2}((-1,0)\times\mathbb{R}^d)$. 
		
		Define $u_2^\varepsilon=u^{(\varepsilon)}-u_1^\varepsilon$. Then $u_2^\varepsilon \in \Sob{1,2}{2}(Q_{3/4})$ satisfies
		\begin{equation*}
			\left\{
			\begin{aligned}
				\partial_t u_2^\varepsilon-D_i(a^{ij}_{(\varepsilon)} D_j u_2^\varepsilon)+\nabla \pi_2^\varepsilon &= \Div(\boldF^{(\varepsilon)}),\\
				\Div u_2^\varepsilon &= g^{(\varepsilon)}
			\end{aligned}
			\right.
		\end{equation*}
		in $Q_{3/4}$. Then following the previous argument as in the nondivergence form case, we can show that $u_2^\varepsilon \in C^0_t C^{1,\alpha}(Q_{1/2})$ for some $\alpha \in (0,1)$. Moreover, it follows from \eqref{eq:boundedness-gradient-estimate} that 
		\begin{align*}
			\norm{Du_2^\varepsilon}{\Leb{\infty}(Q_{1/4})}&\apprle \norm{u_2^\varepsilon}{\Leb{\infty}\Leb{2}(Q_{3/4})}+\norm{Du_2^\varepsilon}{\Leb{2}(Q_{3/4})}+I_{\Rho_f}(1)+I_{\widehat{\Rho}_{g}}(1)\\
			&\apprle \norm{u}{\Leb{\infty}\Leb{2}(Q_1)}+\norm{\boldF}{\Leb{2}(Q_1)}+\norm{g}{\Leb{2}(Q_1)}\\
			&\relphantom{=}+ \norm{u_1^\varepsilon}{\Leb{\infty}\Leb{2}(Q_{3/4})}+\norm{Du_1^\varepsilon}{\Leb{2}(Q_{3/4})}+I_{\Rho_f}(1)+I_{\widehat{\Rho}_{g}}(1).
		\end{align*}
		
		By \eqref{eq:u1-epsilon-div-form} and \eqref{eq:u1-epsilon-Linfty-L2}, we have 
		\[ \norm{u_1^\varepsilon}{\Leb{\infty}\Leb{2}(Q_{3/4})}+\norm{Du_1^\varepsilon}{\Leb{2}(Q_{3/4})}\rightarrow 0\]
		as $\varepsilon\rightarrow 0+$. Hence $\norm{Du_2^\varepsilon}{\Leb{\infty}(Q_{1/4})}$ is bounded by a constant which is independent of $\varepsilon$. On the other hand, it follows from \eqref{eq:u1-epsilon-div-form} that $Du_1^\varepsilon \rightarrow 0$ pointwise a.e. on $Q_{1/4}$. Therefore, we get 
		\[ \norm{Du}{\Leb{\infty}(Q_{1/4})}\apprle \norm{u}{\Leb{\infty}\Leb{2}(Q_1)}+\norm{\boldF}{\Leb{2}(Q_1)}+\norm{g}{\Leb{2}(Q_1)}+I_{\Rho_{\boldF}}(1)+I_{\widehat{\Rho}_{g}}(1).\] 
		The proof is completed.

	\bibliographystyle{amsplain}
	%\bibliography{references}
	\providecommand{\bysame}{\leavevmode\hbox to3em{\hrulefill}\thinspace}
\providecommand{\MR}{\relax\ifhmode\unskip\space\fi MR }
\providecommand{\MRhref}[2]{%
  \href{http://www.ams.org/mathscinet-getitem?mr=#1}{#2}
}
\providecommand{\href}[2]{#2}

\end{document}